\documentclass[10pt]{amsart}
\usepackage{amssymb,amsmath}
\newtheorem{thm}{Theorem}

\newtheorem{defn}[thm]{Definition}
\newtheorem{lemma}[thm]{Lemma}
\newtheorem{cor}[thm]{Corollary}
\newtheorem{note}[thm]{Note}

\numberwithin{equation}{section}

\newcommand{\al}{\alpha}

\begin{document}
\title[Hypoellipticity for systems in Carnot Groups]{Hypoellipticity for linear degenerate elliptic systems in Carnot Groups and
applications}

\begin{abstract}
We prove that if $u$ is a weak solution to a constant coefficient
system (with strong ellipticity assumed along the horizontal
direction) in a Carnot group (no restriction on the step), then
$u$ is actually smooth.  We then use this result to develop
blow-up analysis to prove a partial regularity result for weak
solutions of certain non-linear systems.
\end{abstract}

\author{Emily Shores}
\address{University of Arkansas\\
Fayetteville, AR 72701, USA } \email[Emily
Shores]{eshores23@yahoo.com}

\maketitle

\section{Introduction}
Carnot groups are relatively simple models of sub-Reimannian
manifolds.  In recent years there has been intensive study of the
regularity theory for weak solutions of non-linear degenerate
elliptic systems in this setting. This theory relies on the
hypoellipticity for the corresponding constant coefficient system;
the main purpose of this paper is to prove a regularity result for
weak solutions of the constant coefficient system
\begin{equation}\label{1.1}
\sum_{ \beta = 1}^N \sum _{i,j=1}^m X_i (A_{ij}^{ \alpha \beta}
X_j u^ {\beta} + f_i^{ \alpha}) = f^{ \alpha}\, in \Omega \subset
G\,,
\end{equation}
where $\Omega$ is an open set of a Carnot group G and
$A^{\alpha\beta}_{i,j}$ satisfies the coercivity condition
$$\sum_{\alpha,\beta=1}^N \sum_{{i,j}=1}^m
A^{\alpha\beta}_{i,j} \xi^{\alpha}_i\, \xi^{\beta}_j \geq \lambda
|\xi|^2\,, \qquad  \textrm{for all}\,\xi \in \mathbb{R}^{mN}\,,$$
and $f_i^{ \alpha}$, $f^{ \alpha}$ are smooth functions defined in
$\Omega$. Here $X_i, \ldots, X_m$ refer to the horizontal vector
fields (or rather differentiation along the first layer of the Lie
algebra stratification) in G. We consider weak solutions of
$\eqref{1.1}$ in the horizontal Sobolev space
$S^{1,2}_{loc}(\Omega)$ which consists of $L^2_{loc}(\Omega)$
functions having horizontal derivatives of order one in
$L^2_{loc}(\Omega)$. In particular, we prove the following:

\begin{thm}(Main Theorem)\label{MT}
Let $G$ be a Carnot group of step $r$ and $\Omega \subset G$ an
open, bounded set.  If $u \in S^{1,2}_{loc}(\Omega)$ is a weak
solution to $\eqref{1.1}$ and $f^\alpha$, $f_i^\alpha$ are smooth
functions, then for any ball $B$ such that $2B \subset \Omega$ the
following inequality holds:
\begin{multline}\label{1.2}
\nonumber \parallel X^{I_1}X^{I_2} \cdots X^{I_r}u
\parallel_{S^{1,2}(B)} \\ \leq  C \Bigg( \parallel u
\parallel_{S^{1,2}(2B)} + \parallel \tilde{f}
\parallel_{L^2(2B)} + \parallel \tilde{f_i}
\parallel_{L^2(2B)} \Bigg)\,.
\end{multline}

Here, $X^{I_1} \cdots X^{I_r}$ represent differentiation of
indefinite order in each of the layers, $1, \ldots, r$. Also, we
use the notation $\tilde{f}$, $\tilde{f_i}$ to represent high
order derivatives (possibly along every layer) on the original $f$
and $f_i$.
\end{thm}

The classical method to prove higher regularity of solutions
consists in differentiating the system several times and applying
$L^2$ energy estimates at each step. In Carnot groups this
argument is somewhat different. Carnot groups are non-commutative
groups; therefore, each time we differentiate the system we gain
increasingly complex non-homogeneous terms, the commutators, whose
$L^2$ norm we have to control. These new terms will not only
involve differentiation on the original non-homogeneous terms but
also on the solution $u$. However, since the commutators of two
vector fields belong to a higher layer, one might expect the
derivatives to eventually move to the highest layer, $r$, no
matter in which order we differentiate the system. This is not
always the case as the following example shows:

Suppose we are in a Carnot group, step 4.  Let $X_4$ represent
differentiation in the center of the group $V^4$, $X_3$
differentiation in $V^3$, and $X_2$ differentiation in $V^2$.
Assume that we can differentiate indefinitely along $X_4$ and
$X_3$.  This is a reasonable assumption and is easily shown.  We
wish to then try to differentiate our system a total of three
times, once along each of the directions $X_4$, $X_3$, and $X_2$,
and then show that $X_2X_3X_4u^\beta$ still solves our system.

Without including all the details we will explain where the
problem occurs.  As we try to gain an $L^2$ estimate on
$X_2X_3X_4u^\beta$, we will eventually get to an estimate where
the right hand side includes the sum of the $L^2$ norm of the
following terms\footnote{This term is the non-homogeneous term
that appears in place of $f_i$ when we differentiate along $X_3$
and $X_4$.}(along with similar non-homogeneous terms):
$$Z(A^{\alpha\beta}_{i,j}[X_3,X_j]X_4u + X_3X_4f_i)\quad,X_j \in V^1\,,$$
where $Z$ represents differentiation along each of the vector
fields $X_2$, $X_3$, and $X_4$.  However, the only direction that
poses a problem is $Z=X_2$ (as the other directions are covered
under the assumption). We are assuming that $f_i$ is smooth, so
the second part of the sum, $Z (X_3X_4f_i)$, is bounded in
$L^2_{loc}$. We need to prove an $L^2$ estimate on differentiation
along $X_2$ yet we have a term where this derivative appears on
both the right hand and left hand sides. Therefore, assuming
differentiation along the higher layers is not a sufficient
hypothesis.

What we found is that unless we use the "right" algorithm to
differentiate the system, then phenomena like the one illustrated
above may happen. Referring to this example, one sees that in the
first term of the sum differentiation along the vector field $X_3$
is absorbed into the commutator term.  Thus, differentiation has
essentially moved to the next highest layer once the rules of
commutators are applied. Throughout the paper we refer to this
aspect as a shift in the derivatives to the right.  Moreover, one
will notice that since we no longer have differentiation along
$X_3$ then all of the derivatives lie in only the first and last
layers. This is not only true for the specific example above but
it holds even in the general case; we show that as we apply our
algorithm, the order of the derivatives begins to decrease in the
middle layers until eventually all differentiation shifts to the
first and last layers (see Theorem $\ref{T1}$). The remaining term
is then shown to be bounded above by the $L^2$ norm of a term with
less derivatives than what we start with (see Theorem $\ref{T2}$).
Through an iteration argument we show that this is eventually
bounded above by the $L^2$ norm of $u$. The difficulty in devising
the algorithm is that it must work regardless of what layer one is
differentiating along and regardless of the order of the
derivatives.

As an immediate consequence of the main theorem and of Sobolev's
Embedding Theorem, we have the following:

\begin{cor}\label{MC}
Let $G$ be a Carnot group of step $r$ and $\Omega \subset G$ an
open, bounded set.  If $u \in S^{1,2}_{loc}(\Omega)$ is a weak
solution to $\eqref{1.1}$ and $f^\alpha$, $f_i^\alpha$ are smooth
functions, then $u$ is smooth.
\end{cor}

In the case of scalar equations, Corollary $\ref{MC}$ follows from
a celebrated result of H\"ormander.  In 1967, H\"ormander [H]
studied the partial differential operator $P = \sum_{j=1}^r X_j^2
+ X_0 + c$, where $X_0, \ldots, X_r$ are smooth vector fields in
$\mathbb{R}^n$. He proved that if the vector fields and all of
their commutators generate the whole space, then P is
hypoelliptic. This work, along with the papers of J.J. Kohn [K],
Folland-Stein [FS], Folland [F], and Rothschild-Stein [RS], allow
one to prove the $W^{2,2}$-estimates (and thus the
hypoellipticity) of diagonal systems. We also refer the reader to
the recent papers of Xu and Zuily [XZ] where quasilinear
subelliptic systems are studied, and of Jost and Xu [JX] where
subelliptic harmonic maps are studied. Whereas the above results
address the diagonal case, they do not cover the non-diagonal
case. In this regard, following a highly technical argument the
$W^{2,2}$ estimates can be derived from the analysis of
pseudo-differential operators on homogeneous groups developed in
the papers of [CGGP], [T], and [G].  An advantage of the ideas
presented in the present paper is they may be more familiar to
those working in pde's:  We use fractional order difference
quotients in order to establish differentiation once in any
direction.  Moreover, this method can be applied also to
non-linear systems which cannot be reduced to linear systems (see,
e.g., Theorem 3.9 in [CG]), whereas the method using
pseudo-differential operators cannot.

Our main theorem is a generalization of some of the results proved
in [CG].  In particular, we use a similar approach to establish
hypoellipticity: Roughly speaking we first show that we can
differentiate $\eqref{1.1}$ once in any direction, $Z$. The method
of proof is analogous to the one in [C1], so the proof will be
sketched only. In the step 2 case (see [CG]), once it is
established that the system is differentiable once in any
direction, then indefinite differentiation follows immediately by
an iteration argument. This is not the case for Carnot groups of
arbitrary step, and this is where most of the work in the present
paper lies. The majority of the paper will be devoted to proving
the main theorem, which will directly give us that $u$ has bounded
Sobolev norm of any order. Once this is done we apply the Sobolev
Embedding Theorem (see [F]), to conclude that $u$ is smooth.

The main motivation for our main result comes from non-linear
regularity theory.  We can prove a partial regularity result for
weak solutions of the non-linear system

\begin{equation}\label{1.3}
\sum_{\beta=1}^N \sum_{{i,j}=1}^m {X_i({A_{i,j}^{ \alpha,
\beta}}(x,u){X_ju^ \beta}) =0},\qquad \alpha=1, \ldots, N, \qquad
x \in \Omega,
\end{equation}
where  $\Omega \subset G$ is an open set, $u \in
S^{1,2}_{loc}(\Omega,\mathbb{R}^N)$, and ${A_{i,j}^{ \alpha,
\beta}}(x,u)$ are bounded continuous or uniformly continuous
functions satisfying for a.e. $x \in \mathbb{R}^N$, $u \in
\mathbb{R}$
\begin{eqnarray*}
\sum_{ \alpha, \beta=1}^N \sum_{{i,j}=1}^m {A_{i,j}^{ \alpha,
\beta}}(x,u) {\xi}_i^ \alpha {\xi}_j^ \beta \geq \lambda \vert \xi
\vert ^2, \qquad \xi \in {\mathbb{R}}^{mN}.
\end{eqnarray*}
In fact, we have
\begin{cor}\label{Gi}
If $u\in S^{1,2}_{loc}(\Omega,\mathbb{R}^N)$ is a weak solution to
$\eqref{1.3}$ then there exists an open set $\Omega_0 \subset
\Omega$ such that $u$ is H\"older continuous in $\Omega_0$.
Moreover, the Haar measure of $\Omega \setminus \Omega_0$ is zero.
\end{cor}

This corollary extends to the Carnot group setting a celebrated
result of Giusti and Miranda [GM].  Since much of the elliptic and
degenerate elliptic non-linear regularity theory is based on
elliptic linear estimates, we can consider this result as just a
sample of what can actually be proven using the regularity of the
constant coefficient system (see, for instance, [CG] and [Gi]).
The proof relies on a blow up argument (see [Gi]) and is very
similar to the work done for Carnot groups of step $r=2$ in [CG];
we will briefly describe the argument in section 4, and we then
refer the reader to the papers [CG] and [Gi] for further details
of the proof.

\textit{Acknowledgements} The results in this paper are part of
the authors Ph.D Dissertation at the University of Arkansas.  I
would like to thank my advisor, L. Capogna, for suggesting the
hypoellipticity problem and for his advice and encouragement
throughout.

\section{Preliminaries}

A Carnot group of step $r\geq1$ is defined to be a simply
connected Lie Group $G$ with a decomposition of its lie algebra
$g$ as a vector sum $g=V^1 \oplus V^2 \oplus \cdots \oplus V^r$.
This decomposition is called a stratification of $g$ (of length
$r$) if:  $[V^1,V^j]=V^{j+1}$ for $1 \leq j <r$ and $[V^j,V^r]=0$
for $j \geq r$.  The length of the stratification then corresponds
to the step of the group $G$.  In general, let $X_{i,k}$ denote a
left-invariant basis of $V^k$, where $1 \leq k \leq r=$ step of
$G$ and $1 \leq i \leq m_k=$ dimension of $V^k$. For simplicity,
set $X_i=X_{i,1}$ and $m=m_1$. By the horizontal layer we mean all
of the vectors in the first layer $V^1$. Then we can let $X= \{
X_1, \ldots ,X_m \}$ denote a left-invariant basis for $V^1$. For
a function $u=(u^1, \ldots ,u^N):G \rightarrow {\mathbb{R}}^N$ we
set ${\{Xu \}}_{i,j}=X_iu^j$ to denote the Jacobian of $u$ with
respect to the basis $X$.  We will say that ${\{Xu \}}$ is the
horizontal Jacobian of $u$ since it refers to differentiation in
the horizontal direction only. Following from the fact that the
exponential map $\textrm{exp}: g\rightarrow G$ is a global
diffeomorphism, we can use exponential coordinates on G. We say
that $P \in G$ has coordinates $(p_{i,k})$, for $1 \leq k \leq r$
and $1 \leq i \leq m_k$, where $P=\textrm{exp}
\big(\sum_{k=1}^r\sum_{i=1}^{m_k}p_{i,k}X_{i,k}\big)$.

Carnot groups equipped with the Carnot-Caratheodory metric (see,
e.g. [H]) behave like the Euclidean metric since natural dilations
and translations can be defined.  However, Euclidean spaces are
abelian, and Carnot groups are non-abelian in general. In this
setting, the formula for dilations is given by $\delta_s(P) =
\textrm{exp}(\sum_{k=1}^r \sum_{i=1}^{m_k} s^kp_{i,k}X_{i,k})$,
for $s > 0$ and $P \in G$. It is worth noting that Euclidean
spaces are indeed abelian Carnot groups of step $r=1$. The
simplest example of a non-abelian Carnot group is the Heisenberg
group, ${\mathbb{H}}^1$, which is a Carnot group of step $2$ with
dim $V^2 =1$ and dim $V^1 =2$.

Next we define the pseudo-distance and the gauge balls (see [F]).
First, for $P,Q \in g$ we let $|P|^{2r!} = \sum_{k=1}^r
(\sum_{i=1}^{m_k} |p_{i,k}|^2)^{{r!}/k}$.  Then we have $d(P,Q) =
|Q^{-1}P|$.  In general $d$ does not satisfy the triangle
inequality and $d$ is therefore not a metric.  However, we refer
to $d$ as a gauge metric (or distance).  Second, we use the
psuedo-distance defined above to define the gauge balls.  We have
$B(P,r):= \{x \in G|d(P,x) < r\}$.  We also have that $|B(P,r)| =
\omega_GR^Q$, where $\omega_G = |B(e,1)|$, $e$ is the group
identity, and $Q = \sum_{k=1}^r km_k$ is the so-called homogeneous
dimension of $G$ ([F]).

Next, we remind the reader of the definition of horizontal Sobolev
spaces.
\begin{defn}
For $k \in \mathbb{N}$, $1\leq p \leq \infty$, and for $\Omega
\subset G$, we let ${\mathcal{S}^{k,p}_{loc}}(\Omega)$ represent
the set of functions $f: \Omega \rightarrow \mathbb{R}^N$ such
that the components of $f$ are in $L^p_{loc}(\Omega)$ and all of
the horizontal derivatives of the components of $f$ of order up to
$k$ are in $L^p_{loc}(\Omega)$.
\end{defn}

The ${\mathcal{S}^{k,p}_{loc}}(\Omega)$ norm is then given by
$\parallel f \parallel_{{\mathcal{S}^{k,p}}}=\parallel f
\parallel_{L^p} + \sum_{l=1}^k \sum_{I \in(1,\ldots ,m_1)^l}
\parallel X_{i_1,1}X_{i_2,1}\ldots X_{i_l,1}f \parallel_{L^p}$.
Note that if $u \in {\mathcal{S}^{k,p}_{loc}}(\Omega)$ for all $k$
then we also have $u \in {\mathcal{W}^{k,p}_{loc}}(\Omega)$ for
all $k$, where ${\mathcal{W}^{k,p}_{loc}}(\Omega)$ represents the
usual Euclidean Sobolev space.

\begin{defn}
A function $u \in {\mathcal{S}^{1,2}_{loc}}(\Omega)$ is a weak
solution of \eqref{1.1} if we have the following identity for each
$\phi \in C^{^\infty_0}(\Omega)$:
\begin{equation}
\sum_{ \beta = 1}^N \sum _{i,j=1}^m \int_{\Omega} \Big(A_{ij}^{
\alpha \beta} X_j u^ {\beta} + f_i^{ \alpha}\Big)(p) X_i
\phi^{\alpha}(p)\,dp = \int_{\Omega} f^{
\alpha}\phi^{\alpha}(p)\,dp \,.
\end{equation}
\end{defn}

Before we prove the main theorem, we recall the following results.
These show that we can differentiate our system \eqref{1.1} once
in any direction, $Z$, and that $Zu$ is still a solution to the
system.

\begin{thm}\label{A}
Let G be a Carnot group of step r and $\Omega \subset G$ an open
set, and $u \in S^{1,2}_{loc}(\Omega)$ a weak solution to
(\ref{1.1}).  If $f^\alpha$, $X_{r,i_0}f^\alpha$, $f_i^\alpha$,
$X_{r,i_0}f_i^\alpha \in L^2_{loc}(\Omega)$ for every $1 \leq i_0
\leq m_r$ and $\alpha = 1, \ldots, N$, then for every $X_{r,i_0}
\in V^r$ one has $$X_{r,i_0}u \in S^{1,2}_{loc}(\Omega)\,.$$
Furthermore, for every pseudo-ball $B(p_o,2R) \subset \Omega$, the
following estimate holds for $X_{r,i_0}u$

\begin{eqnarray}\
 \nonumber \parallel X_{r,i_0}u \parallel_{S^{1,2}_{loc}(B(p_0,R))}
 &\leq&
 C \Bigg( \parallel u
\parallel_{S^{1,2}_{loc}(B(p_0,2R))} + \parallel f
\parallel_{L^2(B(p_0,2R))} + \parallel X_{r,i_0}f
\parallel_{L^2(B(p_0,2R))}\\ &+& \sum_{i=1}^m \Big[ \parallel
f_i\parallel_{L^2(B(p_0,2R))} + \parallel
X_{r,i_0}f_i\parallel_{L^2(B(p_0,2R))}\Big] \Bigg)
\end{eqnarray}
for some constant $C >0$ depending only on $g$ and the coercivity
condition.  Moreover, $X_{r,i_0}u$ is a weak solution to the
system
\begin{equation}
\sum_{\beta =1}^N \sum_{i,j=1}^m X_i \big(A_{ij}^{\alpha
\beta}X_j(X_{r,i_0}u)^\beta + X_{r,i_0}f_i^\alpha\big) =
X_{r,i_0}f^\alpha\,,
\end{equation}
for every $\alpha = 1, \ldots, N$.
\end{thm}

\begin{thm}\label{B}
Let $u \in S^{1,2}_{loc}(\Omega)$ be a weak solution to
(\ref{1.1}), $\Omega \subset G$ be an open set, and $G$ a Carnot
group, step $r$.  Assume
$\omega_{\tilde{k}}=X_{\tilde{k},l}u^\alpha \in
S^{1,2}_{loc}(\Omega)$ for every $1\leq k_0<\tilde{k}$, and $1\leq
l \leq m_{\tilde{k}}$ such that $\omega_{\tilde{k}}$ satisfies

\begin{eqnarray}\label{B1}
\parallel \omega_{\tilde{k}}\parallel_{S^{1,2}_{loc}(B(p_0,R))}
&\leq& C \Bigg[ \parallel u \parallel_{S^{1,2}_{loc}(B(p_0,2R))}
\nonumber  +
\parallel f \parallel_{L^2(B(p_0,2R))} \\&+&
\sum_{j=1}^{m_{\tilde{k}}}\sum_{k=\tilde{k}}^r \parallel
X_{k,j}f\parallel_{L^2(B(p_0,2R))} + \sum_{i=1}^m
\parallel f_i
\parallel_{L^2(B(p_0,2R))}
\\&+& \nonumber \sum_{j=1}^{m_{\tilde{k}}}\sum_{k=\tilde{k}}^r
\sum_{i=1}^m
\parallel X_{k,j}f_i \parallel_{L^2(B(p_0,2R))}\Bigg]
\end{eqnarray}
for any ball $B(p_0,2R)\subset \Omega$.  Further, for $\alpha = 1,
\ldots, N$, if $f^\alpha$, $f_i^\alpha$, $X_{k_0,i_0}f^\alpha$,
$X_{k_0,i_0}f_i^\alpha \in L^2_{loc}(\Omega)$ then we have
$\omega_{k_0} = X_{k_o,i_o}u \in S^{1,2}_{loc}(\Omega)$.

Moreover, $\omega_{k_0}$ is a weak solution to

\begin{eqnarray}\label{B2}
&&\sum_{\beta=1}^N\sum_{i,j=1}^m X_i \Bigg( A_{ij}^{\alpha \beta}
X_j X_{k_0,i_0}u^\beta \,+\, A_{ij}^{\alpha \beta} [
X_{k_0,i_0},X_j] u^\beta \,+\, X_{k_0,i_0}f_i^\alpha \Bigg) \\ &=&
\nonumber  X_{k_0,i_0}f^\alpha \,+\,
\sum_{\beta=1}^N\sum_{i,j=1}^m [X_i,X_{k_0,i_0}]
\big(A_{ij}^{\alpha \beta} X_j u^\beta \,+\, f_i^\alpha \big)
\end{eqnarray}
for $\alpha = 1, \ldots, N$\,, and $\eqref{B1}$ holds for
$\tilde{k}=k_0$.
\end{thm}

\begin{proof}[Sketch of proof]

The method of proof in [C1] still holds in this setting of
systems, and we therefore refer the reader to that paper for the
details of the proofs. We will need estimates on the Lebesgue norm
of fractional derivatives of functions in the direction of
commutators, so we must first introduce the following notation.
Let $\Omega$ be an open subset of G, $Z \in g$, $\omega \in
L^2(\Omega)$ with compact support in $\Omega$, and $\alpha \in
(0,1)$.  We define the seminorm $$|\omega|^2_{Z,\alpha} =
\sup_{|h| <\epsilon_0} \int_\Omega |h|^{-2\alpha}
|\omega(ze^{hZ})-\omega(z)|^2\,dz\,,$$ where $\epsilon_0$ is
chosen sufficiently small.  Then we can express the $L^2$-norm of
the fractional derivative of $\omega$ along the direction
$\partial_{p_j,l}$ in terms of exponential coordinates by the
formula
\begin{eqnarray*}
\parallel \partial_{p_j,l}^\alpha \omega \parallel_{L^2(G)} =
&\int_G&
|h|^{2\alpha}|\hat{\omega}(p_{1,1},\ldots,p_{j-1,l},h,p_{j+1,l},\ldots,p_{m_r,r})|^2\\&&
dp_{1,1} \ldots dp_{j-1,l} dh dp_{j+1,l} \ldots dp_{m_r,r}\,,
\end{eqnarray*}
where we have denoted by $\hat{\omega}$ the partial Fourier
transform in the variable $p_{j,l}$. Next, we use the following
theorems of Peetre [P] and H\"ormander (Theorem 4.3 in [Ho]),
along with the Energy inequality:

\begin{thm}[Peetre]\label{Peetre}
Let $G$ be a Carnot group of step $r$, let $0 < \beta < \alpha <
1$, and $ \omega \in C^{\infty}_0(g)$. Then there exists positive
constant $C = C ( \alpha, \beta, N)$ such that
\begin{eqnarray*}
C \parallel \partial^{\beta}_{p_j,l} \omega \parallel_{L^2(G)}
\leq | \omega|_{ \partial_{p_j,l, \alpha}} \leq C^{-1} \parallel
\partial^{\alpha}_{p_j,l} \omega \parallel_{L^2(G)}
\end{eqnarray*}
where $p \in G$ has the coordinates $p_{i,k}$, for $1 \leq i \leq
m_k$, $1 \leq k \leq r$.
\end{thm}

\begin{thm}[H\"ormander]\label{Hormander}
Let $\omega \in C^{\infty}_0(G)$.  For $ 1 \leq k \leq r$ , $1
\leq i \leq m_k$ one has
\begin{eqnarray*}
| \omega|_{X_{i,k, \frac{1}{k}}} \leq C \sum_{j=1}^m |
\omega|_{X_j,1} + \parallel \omega \parallel_{L^2(G)},
\end{eqnarray*}
for some positive constant $C$, and for $G$, a Carnot group of
step $r$.
\end{thm}

\begin{lemma}[Energy inequality]
Let $G$ be a Carnot group of step $r$, $\Omega \subset G$ be an
open set (bounded).  If $u \in S^{1,2}_{loc} (\Omega)$ is a weak
solution to (\ref{1.1}) in $\Omega$, with the assumption that
there exists $\lambda > 0$ such that for every $x \in \Omega$, one
has
$$\sum_{\alpha, \beta = 1}^N \sum _{i,j=1}^m A^{\alpha
\beta}_{i,j}(x,u)\, \xi_i^\alpha \xi_j^\beta \geq \lambda
(M)|\xi|^2,$$ then the following Caccioppoli-type inequality holds
for $2B= B(p_0,2r)\subset \Omega:$
\begin{equation}
\int_{B(p_0,r)}|Xu|^2\,dp \leq \frac{C}{r^2} \int_{2B}|u|^2\,dp +
C \int_{2B}\Big(|f|^2 + \sum_{i=1}^m|f_i|^2\Big)\,dp
\end{equation}
\end{lemma}

Roughly speaking, to prove Theorem $\ref{A}$ we consider
fractional difference quotients of $u$ in the direction $Z$ of
order $\alpha \in (0,1]$ and apply Theorems $\ref{Peetre}$ and
$\ref{Hormander}$ in order to show that we can actually consider
difference quotients of order $1$. We express the fractional
difference quotient using the formula $u_{(Z,\alpha)}(p)=
\frac{u(pe^{sZ})-u(p)}{|s|^\alpha}$. Utilizing the Caccioppoli
inequality gives us the result. Theorem \ref{B} follows from
Theorem \ref{A} by using an iteration argument (based on the layer
being differentiated in) to give us differentiation of order one
in any direction.
\end{proof}

\section{Proof of Main Theorem}

In this section we will show that if $u \in S^{1,2}_{loc}(\Omega)$
is a weak solution to the constant coefficient system \eqref{1.1}
then $u$ is smooth.

\subsection{\textbf{Notation}}
For every $h_k \in \mathbb{N}$, and for the multi-indices
$$I_{h_k} = (i_1,i_2, \ldots, i_{h_k}) \in \{1,2, \ldots
m_k\}^{h_k}\,,$$ we define the following terms. Throughout the
paper we let $l-1$ represent the lowest layer that we are
differentiating with respect to.  So to represent differentiating
${(h_{l-1})+1}$ times with respect to this layer only, we set
$$X^{I_{h_{l-1}}}=X_{l-1,i_{(h_{l-1})+1}}X_{l-1,i_{h_{l-1}}}
\cdots X_{l-1,i_1} \quad.$$ Then for each $k > l-1$, i.e. each
layer above $l-1$, the following definition represents taking
$h_k$ derivatives within each of the $k$ layers:
$$X^{I_{h_k}}=X_{k,i_{h_k}}X_{k,i_{(h_k)-1}} \cdots X_{k,i_1} \quad.$$
We have set different notation for the lowest layer, $k=l-1$, that
we are differentiating with respect to; this is simply a matter of
convenience in order to make the computations more clear.

Whereas the above two definitions represent taking multiple
derivatives in one layer at a time, the next two definitions give
us notation to represent taking multiple derivatives in multiple
layers.  For $k \geq l-1$ we set

$$V(k)\,=\,X^{I_{h_k}}X^{I_{h_{k+1}}}\cdots X^{I_{h_r}}u^\beta$$
$$V(k+1)^{k,i_{(h_k)-s}} \,= \,X^{I_{(h_k)-s}}X^{I_{h_{k+1}}}\cdots
X^{I_{h_{r}}}u^\beta \,=\,X^{I_{(h_k)-s}}V(k)\,,\quad \textrm{for
s} \leq\, h_k\,. $$ The difference in the above two definitions is
the second one keeps count of the derivatives in the lowest layer.

Each time we differentiate the system, we end up with
non-homogeneous terms due to the non-commutativity of the group
structure.  After only a few steps into this differentiation
process, one can see that these terms are complicated and quickly
become difficult to work with; we will use the following notation
to define such non-homogeneous terms, with each one being defined
in terms of the previous one in order to simplify the
computations. For every $k \geq l-1$, set

\begin{eqnarray*}
f_i^\alpha(k)^{k-1,i_{h_{k-1}}} &=&  A_{ij}^{ \alpha
\beta}[X_{k-1,i_{h_{k-1}}},X_j]V(k)^{k-1,i_{(h_{k-1})-1}} \\&+&
X_{k-1,i_{h_{k-1}}}f_i^\alpha(k)^{k-1,i_{(h_{k-1})-1}}\,,\\
f^\alpha(k)^{k-1,i_{h_{k-1}}} &=&
X_{k-1,i_{h_{k-1}}}f^\alpha(k)^{k-1,{i_{(h_{k-1})-1}}}
\\ &+&
[X_i,X_{k-1,i_{h_{k-1}}}] \big(A_{ij}^{ \alpha \beta} X_j
V(k)^{k-1,i_{(h_{k-1})-1}}  +
f_i^\alpha(k)^{k-1,i_{(h_{k-1})-1}}\big)\,,\\
f_i^\alpha(k)^{k-1,0}  &=& f_i^\alpha(k+1)^{k,i_{h_k}}\,,\\
f^\alpha(k)^{k-1,0}  &=& f^\alpha(k+1)^{k,i_{h_k}}\,,\\
f_i^\alpha(r)^{r-1,0}  &=&  (f_i^\alpha)^{r,i_{h_r}}  =
X^{I_{h_r}}f_i^\alpha \,,\\
f^\alpha(r)^{r-1,0} &=& (f^\alpha)^{r,i_{h_r}}  =
X^{I_{h_r}}f^\alpha \,.
\end{eqnarray*}

Last, we need a way to represent taking a different number of
derivatives in any one layer than what we started with.  We
introduce the following notation that will represent taking $b$
derivatives within the single layer $V^k$, for $b \leq h_k$ and
$l-1 \leq k \leq r$. We let
$$X^{k,b}=X_{k,i_b}X_{k,i_{(b-1)}} \ldots X_{i_1}\,.$$
\begin{note}To simplify somewhat the heavy notation, we will set  $*
=i_{h_{l-1}}$ throughout the paper.\end{note}

\subsection{\textbf{Results}}
We are assuming that $u\equiv(u^1, \cdots, u^N):\Omega \rightarrow
\mathbb{R}^N \in S^{1,2}_{loc}(\Omega)$ is a weak solution to
\eqref{1.1} for every $\beta = 1, \cdots, N$, and our first aim is
to show that $X^{I_{h_{l-1}}}X^{I_{h_l}} \cdots X^{I_{h_r}}u^\beta
\in W^{1,2}_{loc}(\Omega)$ for every $1 \leq l-1 \leq r$, for
every $\beta$, and that this is a weak solution to the system
$$\sum_{ \beta = 1}^N \sum _{i,j=1}^m X_i \Big(A_{ij}^{ \alpha \beta}
X_j {V(l-1)} + f_i^\alpha(l)^{l-1,i_{(h_{l-1})+1}} \Big)
=f^\alpha(l)^{l-1,i_{(h_{l-1})+1}} \qquad.$$

In order to achieve this goal, our proof is divided into two main
steps which are detailed in the following two theorems.

\begin{thm}\label{T1}
Let $u \in S^{1,2}_{loc}(\Omega)$ be a weak solution to the system
\eqref{1.1} with $f^\alpha, f_i^\alpha \in C^\infty(\Omega)$. If
$\tilde{f}$, $\tilde{f_i}$ are as in the statement of Theorem
\ref{MT}, then we have:
\begin{eqnarray}\label{T1eq1}
\nonumber \parallel X^{I_{h_{l-1}}}X^{I_{h_{l}}} &\cdots&
X^{I_{h_r}}u^\beta
\parallel_{S^{1,2}_{loc}(B)}
\\& \leq & \parallel
X^{J_{h_{l-1}}}X^{J_{h_l}} \cdots X^{J_{h_r}}u^\beta
\parallel_{S^{1,2}_{loc}(2B)} \\
&+& \nonumber \parallel \tilde{f^\alpha} \parallel_{{L^2}(2B)} +
\parallel \tilde{f_i^\alpha}
\parallel_{{L^2}(2B)}\,,
\end{eqnarray}
where either
$$\sum_{i=l-1}^r |J_{h_i}| < \sum_{i=l-1}^r
|I_{h_i}|$$ or $$\sum_{i=l-1}^r |J_{h_i}| = \sum_{i=l-1}^r
|I_{h_i}|\,,$$ and one of the following two things occur:
\begin{itemize}
\item [(i.)] $|J_{h_{l-1}}| < |I_{h_{l-1}}|$ and consequently there exists at least
one $\beta > l-1$ with $|J_{h_\beta}| > |I_{h_\beta}|$
\item [(ii.)] $|J_{h_{l-1}}| = |I_{h_{l-1}}|$ and there exists $\beta >
l-1$ with $|J_{h_\beta}| < |I_{h_\beta}|$ and $|J_{h_{\beta+1}}| >
|I_{h_{\beta+1}}|\,.$
\end{itemize}
\end{thm}

\begin{thm}\label{T2}
Let $u \in S^{1,2}_{loc}(\Omega)$ be a weak solution to the system
\eqref{1.1} with $f ^\alpha$, $f_i^ \alpha \in C^\infty(\Omega)$.
If $\tilde{f}$, $\tilde{f_i}$ are as in the statement of Theorem
\ref{MT}, then we have:

\begin{eqnarray}\label{T2eq1}
\nonumber \parallel X^{I_{h_{l-1}}}X^{I_{h_r}}u^\beta
\parallel_{S^{1,2}_{loc}(B)} & \leq & \parallel
X^{J_{h_{l-1}}}X^{J_{h_l}} \cdots X^{J_{h_r}}u^\beta
\parallel_{S^{1,2}_{loc}(2B)} \\
&+& \parallel \tilde{f^\alpha} \parallel_{{L^2}(2B)} + \parallel
\tilde{f_i^\alpha}
\parallel_{{L^2}(2B)}\,,
\end{eqnarray}
where $$\sum_{i=l-1}^r |J_{h_i}| \leq \sum_{i=l-1}^r
|I_{h_i}|\,,$$ with $|J_{h_{l-1}}|$ always being at least one less
than $|I_{h_{l-1}}|$, and $|J_{h_k}| \geq |I_{h_k}|$ for every
other $k$.
\end{thm}

The crucial step to proving Theorems \ref{T1} and \ref{T2} is the
estimates on the $L^2$ norm of the terms $f^\alpha(l)^{l-1,*}$ and
$f_i^\alpha(l)^{l-1,*}$. We start by showing that when we apply
the definitions of $f^\alpha(l)^{l-1,*}$ and
$f_i^\alpha(l)^{l-1,*}$ we have done one of two things.  Either we
have lessened the number of derivatives in the lowest layer (and
thus lessened the total number of derivatives) or we have kept the
same number in the lowest layer yet shifted the derivatives on
$u^\beta$ somewhere to the right of the lowest layer, without
adding to the total number of derivatives. Iterating Theorem
\ref{T1} will eventually shift all of the derivatives to the first
and last layers, l-1 and r, respectively.  Theorem \ref{T2} tells
us that derivatives in the first and last layer can be bounded
above by an $L^2$ estimate in which we have lessened the number of
derivatives in the lowest layer yet added some to the higher
layers. Iterating these two theorems will then give us that the
$S^{1,2}_{loc}$ norm of $X^{I_{h_{l-1}}}X^{I_{h_l}} \cdots
X^{I_{h_r}}u^\beta$ is bounded above.

\subsection{Auxillary Lemmas}
Recall the definitions of $f^\alpha(k)^{k-1,i_{h_{k-1}}}$ and
$f_i^\alpha(k)^{k-1,i_{h_{k-1}}}$ above.  The following two lemmas
provide estimates on the $L^2$ norms of
$f^\alpha(k)^{k-1,i_{h_{k-1}}}$ and
$f_i^\alpha(k)^{k-1,i_{h_{k-1}}}$ in terms of the $L^2$ norms of
our original $f^\alpha$, ${f_i}^\alpha$, and $u^\beta$.

\begin{lemma}\label{L1}
If $u^\beta$, $f^\alpha$, and $f_i^\alpha$ are the same as in
Theorem \ref{MT}, then we have the following:
\begin{eqnarray*}
\parallel f^\alpha(l)^{l-1,*}\parallel_{L^2}& \leq &
C \Bigg (\parallel X^{l-1,*}X^{I_{h_l}} \cdots X^{I_{h_r}}f^\alpha
\parallel_{L^2}\\& +& \sum_{q+k=*-1}
\parallel  X^{l-1,q}X^{l,1}X_jV(l)^{l-1,k}
\parallel_{L^2}\\
&+& \sum_{s=l+1}^r \parallel W(s) \parallel_{L^2}\\ &+&
\sum_{q+k=*-1}
\parallel X^{l-1,q}X^{l,1} f_i^\alpha(l)^{l-1,k}
\parallel_{L^2}\Bigg)\,,
\end{eqnarray*}

where
\begin{eqnarray*}
W(s) &=& \sum_{q+k=i_{(h_{s-1})-1}} X^{l-1,*}\cdots
X^{s-1,q}X^{s,1}X_jX^{s-1,k}X^{I_{h_s}} \cdots X^{I_{h_r}}u^\beta  \\
&+& \sum_{q+k=i_{(h_{s-1})-1}}
 X^{l-1,*}\cdots X^{s-1,q}X^{s,1} f_i^\alpha(s)^{s-1,k}
\end{eqnarray*}
and the constant C depends on the coefficients.
\end{lemma}
\begin{proof}

Referring to the definition of $f^\alpha(l)^{l-1,*}$ and using the
notation\footnote{Suppose we have $[X_{l,a},X_{m,b}]$.  This is a
linear combination of the vector field that one obtains by adding
the subscripts together, i.e. $X_{(l+m,\cdot)}$.  For simplicity,
whenever we commute two vector fields, say $X_{l,a}$ and
$X_{m,b}$, we will call the new term $X^{l+m,1}$ and drop the
second subscript.} $[X_i,X_{k,i_{h_k}}]= X^{k+1,1}$, we have:
\begin{eqnarray}\label{l1}
\nonumber f^\alpha(l)^{l-1,*} &=& X_{l-1,*}f^\alpha(l)^{l-1,*-1}
\\ \nonumber &+& [X_i,X_{l-1,*}] \Big(A_{ij}^{ \alpha
\beta}X_jV(l)^{l-1,*-1}+ f_i^\alpha(l)^{l-1,*-1} \Big)\\
&=& X_{l-1,*}f^\alpha(l)^{l-1,*-1} + A_{ij}^{ \alpha
\beta}X^{l,1}X_jV(l)^{l-1,*-1} \\ \nonumber &+& X^{l,1}
f_i^\alpha(l)^{l-1,*-1}\quad.
\end{eqnarray}

Next, we rewrite $f^\alpha(l)^{l-1,*-1}$ using this new
representation.  Doing so and then substituting the result back
into (\ref{l1}) yields:
\begin{eqnarray*}
f^\alpha(l)^{l-1,*} &=&
X_{l-1,*}\Big(X_{l-1,*-1}f^\alpha(l)^{l-1,*-2} \\ &+& A_{ij}^{
\alpha \beta}X^{l,1}X_jV(l)^{l-1,*-2} + X^{l,1}
f_i^\alpha(l)^{l-1,*-2}\Big)\\&+& A_{ij}^{ \alpha
\beta}X^{l,1}X_jV(l)^{l-1,*-1} + X^{l,1}
f_i^\alpha(l)^{l-1,*-1}\quad.
\end{eqnarray*}

Iterating this process ($*-2$) more times admits the following:

\begin{eqnarray}\label{2}
f^\alpha(l)^{l-1,*} &=& X^{l-1,*}f^\alpha(l)^{l-1,0} +
\sum_{q+k=*-1}A_{ij}^{ \alpha \beta}
X^{l-1,q}X^{l,1}X_jV(l)^{l-1,k}\nonumber \\
&+& \sum_{q+k=*-1}X^{l-1,q}X^{l,1}f_i^\alpha(l)^{l-1,k} \nonumber \\
&=& X^{l-1,*}f^\alpha(l+1)^{l,i_{h_l}} + \sum_{q+k=*-1}A_{ij}^{
\alpha \beta} X^{l-1,q}X^{l,1}X_jV(l)^{l-1,k}
\nonumber \\
&+& \sum_{q+k=*-1}X^{l-1,q}X^{l,1}f_i^\alpha(l)^{l-1,k}\,.
\end{eqnarray}

So far we have rewritten $f^\alpha(l)^{l-1,*}$ in terms of the
next higher step with our eventual goal being to rewrite it based
on differentiation along the original $f^\alpha$. Once we finish
rewriting this term there will still be other terms that appear,
in particular ones similar to the $f_i^\alpha(l)^{l-1,k}$ above.
However, our next lemma will concern terms of this type, so for
now we will leave these as is.  What we do next is continue the
process above by using (\ref{2}) on each new
$f^\alpha(k)^{k-1,i_{h_{k-1}}}$, for $l+1 \leq k\leq r$.

First adapt equation (\ref{2}) to the term
$f^\alpha(l+1)^{l,i_{h_l}}$ and then rewrite (\ref{2}) using this.
We repeat this process next for $f^\alpha(l+2)^{l+1,i_{h_{l+1}}}$,
and we continue doing this until we last apply it to
$f^\alpha(r)^{r-1,i_{h_{r-1}}}$. Proceeding in this way, one has
the following representation for $f^\alpha(l)^{l-1,*}$:

\begin{eqnarray*}
f^\alpha(l)^{l-1,*} &=& X^{l-1,*}X^{l,i_{h_{l}}} \cdots
X^{r,i_{h_r}}f^\alpha \\ &+& \sum_{s=l+1}^r
\Big(\sum_{q+k=i_{(h_{s-1})-1}}A_{ij}^{ \alpha \beta} X^{l-1,*}
\cdots X^{s-1,q}X^{s,1}X_jV(s)^{s-1,k}\Big)
\\ &+&
\sum_{s=l+1}^r \Big(\sum_{q+k=i_{(h_{s-1})-1}}X^{l-1,*} \cdots
X^{s-1,q}X^{s,1}f_i^\alpha(s)^{s-1,k}\Big) \\ &+&
\sum_{q+k=*-1}A_{ij}^{ \alpha \beta}
X^{l-1,q}X^{l,1}X_jV(l)^{l-1,k}
\\ &+& \sum_{q+k=*-1}X^{l-1,q}X^{l,1}f_i^\alpha(l)^{l-1,k} \\ &=&
X^{l-1,*}X^{l,i_{h_{l}}} \cdots X^{r,i_{h_r}}f^\alpha +
\sum_{s=l+1}^r W(s) \\ &+& \sum_{q+k=*-1}A_{ij}^{ \alpha \beta}
X^{l-1,q}X^{l,1}X_jV(l)^{l-1,k} \\ &+&
\sum_{q+k=*-1}X^{l-1,q}X^{l,1}f_i^\alpha(l)^{l-1,k}\,.
\end{eqnarray*}

Lastly, take the $L^2$ norm of both sides to obtain the desired
result.
\end{proof}

\begin{lemma}\label{L2}
If $u^\beta$, $f^\alpha$, and $f_i^\alpha$ are the same as in
Theorem \ref{MT}, then we have the following:
\begin{eqnarray*}
\parallel f_i(l)^{l-1,*}\parallel_{L^2}& \leq &
C \Bigg(\parallel X^{l-1,*}X^{l,i_{h_l}} \cdots X^{r,i_{h_r}}f_i
\parallel_{L^2} + \sum_{s=l+1}^r \parallel T(s) \parallel_{L^2}\\
&+& \sum_{q+k=*-1} \parallel X^{l-1,q}X_lV(l)^{l-1,k}
\parallel_{L^2}\Bigg)\,,
\end{eqnarray*}
where
$$T(s) = \sum_{q+k=i_{(h_{s-1})-1}}
X^{l-1,*} \cdots X^{s-1,q}X_sV(s)^{s-1,k}\,.$$
\end{lemma}
\begin{proof}
Referring to the definition of $f_i(l)^{l-1,*}$ and using the
notation $[X_i,X_{k,i_{h_k}}]= X^{k+1,1}$, we have:
\begin{eqnarray*}
f_i(l)^{l-1,*}  &=& A_{ij}^{ \alpha
\beta}[X_{l-1,*},X_j]V(l)^{l-1,*-1} + X_{l-1,*}f_i(l)^{l-1,*-1}\\
&=& A_{ij}^{ \alpha \beta} X^{l,1} V(l)^{l-1,*-1} +
X_{l-1,*}f_i(l)^{l-1,*-1}\,.
\end{eqnarray*}

Next we apply this representation to $f_i(l)^{l-1,*-1}$ and
substitute the result into the equality above to obtain:

\begin{eqnarray*}
f_i(l)^{l-1,*}  &=&  A_{ij}^{ \alpha \beta} X^{l,1} V(l)^{l-1,*-1}
+ X_{l-1,*}\Big(A_{ij}^{ \alpha \beta} X^{l,1} V(l)^{l-1,*-2}
\\ &+& X_{l-1,*-1}f_i(l)^{l-1,*-2}\Big) \\ &=&
A_{ij}^{ \alpha \beta}X^{l,1}V(l)^{l-1,*-1} +A_{ij}^{
\alpha \beta} X_{l-1,*}X^{l,1} V(l)^{l-1,*-2} \\
&+& X_{l-1,*}X_{l-1,*-1}f_i(l)^{l-1,*-2}\,.
\end{eqnarray*}

Iterating this $(*-2)$ more times we get:
\begin{eqnarray*}
f_i(l)^{l-1,*}  &=& \sum_{q+k=*-1}A_{ij}^{ \alpha
\beta}X^{l-1,q}X^{l,1}V(l)^{l-1,k} + X^{l-1,*}f_i(l)^{l-1,0}
\\ &=& \sum_{q+k=*-1}A_{ij}^{ \alpha
\beta}X^{l-1,q}X^{l,1}V(l)^{l-1,k} +
X^{l-1,*}f_i(l+1)^{l,i_{h_l}}\,.
\end{eqnarray*}

Continue by first applying the argument above to
$f_i(l+1)^{l,i_{h_l}}$, then to $f_i(l+2)^{l+1,i_{h_{l+1}}}$, etc,
and last to $f_i(r)^{r-1,i_{h_{r-1}}}$.  Proceeding in this way,
one obtains the following representation for $f_i(l)^{l-1,*}$:
\begin{eqnarray*}
f_i(l)^{l-1,*} &=& \sum_{q+k=*-1} A_{ij}^{ \alpha \beta}
X^{l-1,q}X^{l,1}V(l)^{l-1,k} \\ &+& \sum _{s=l+1}^r
\Big(\sum_{q+k=i_{(h_{s-1})-1}}A_{ij}^{ \alpha \beta} X^{l-1,*}
\cdots X^{s-1,q}X^{s,1}V(s)^{s-1,k}\Big) \\ &+&
X^{l-1,*}X^{l,i_{h_l}} \cdots X^{r,i_{h_r}}f_i \\ &=&
\sum_{q+k=*-1} A_{ij}^{ \alpha \beta} X^{l-1,q}X^{l,1}V(l)^{l-1,k}
+ \sum _{s=l+1}^r T(s) \\ &+& X^{l-1,*}X^{l,i_{h_l}} \cdots
X^{r,i_{h_r}}f_i\,.
\end{eqnarray*}
Last, take the $L^2$ norms of both sides to complete the proof.
\end{proof}

We can apply lemma \ref{L2} to rewrite lemma \ref{L1}.  As a
direct consequence we have:

\begin{lemma}\label{L3}
If $u^\beta$, $f^\alpha$, and $f_i^\alpha$ are the same as in
Theorem \ref{MT}, then we have the following:
\begin{eqnarray*}
\parallel  &f(l) ^{l-1,*}& \parallel_{L^2}\\ & \leq &
C \Bigg(\sum_{q+k=*-1}
\parallel X^{l-1,q}X^{l,1}X_jV(l)^{l-1,k}
\parallel_{L^2} \\ &+& \sum_{s=l+1}^r \Big( \sum_{q+k=i_{(h_{s-1})-1}}
\parallel X^{l-1,*}\cdots X^{s-1,q}X^{s,1}X_jX^{s-1,k}X^{s,i_{h_s}} \cdots
\\ & \cdots & X^{r,i_{h_r}}u\parallel_{L^2}\Big) \\ &+&
\sum_{q+y+z=*-2}
\parallel X^{l-1,q}X^{l,1}
X^{l-1,y}X^{l,1}V(l)^{l-1,z}
\parallel_{L^2}\\ &+& \sum_{q+k=*-1}\Big(\sum_{y+z=i_{(h_{s-1})-1}}\Big(
\sum_{s=l+1}^r
\parallel X^{l-1,q}X^{l,1} X^{l-1,k}X^{l,i_{h_l}} \cdots \\ &\cdots&
X^{s-1,y}X^{s,1}V(s)^{s-1,z}\parallel_{L^2}\Big)\Big)\\&+&
\sum_{q+k=i_{(h_{s-1})-1}}\Big(\sum_{y+z=k-1} \Big( \sum_{s=l+1}^r
\parallel X^{l-1,*}\cdots
X^{s-1,q}X^{s,1}X^{s-1,y} X^{s,1}V(s)^{s-1,z}\parallel_{L^2}\Big)\Big) \\
&+&
\sum_{q+k=i_{(h_{s-1})-1}}\Big(\sum_{y+z=i_{(h_{p-1})-1}}\Big(\sum_{s=l+1}^r
\Big( \sum_{p=s+1}^r
\parallel X^{l-1,*} \cdots \\ & \cdots &
X^{s-1,q}X^{s,1}X^{s-1,k}X^{s,i_{h_{s}}} \cdots
X^{p-1,y}X^{p,1}X^{p-1,z}X^{p,i_{h_{p}}}V(p+1)
\parallel_{L^2}\Big)\Big)\Big)
 \\ &+& C(f) + C(f_i)\Bigg)\quad,
\end{eqnarray*}
where $C(f)$ and $C(f_i)$ are terms corresponding to
differentiation on the original $f$ and $f_i$, respectively, and
$C$ is a constant depending on the coefficients.
\end{lemma}

\begin{proof}

The two terms in lemma \ref{L1} to focus on are $X^{l-1,q}X^{l,1}
f_i(l)^{l-1,k}$ and $X^{l-1,*}\cdots X^{s-1,q}X^{s,1}
f_i(s)^{s-1,k}$, so begin by applying lemma \ref{L2} to the first
term:

\begin{eqnarray*}
\sum_{q+k=*-1} & \parallel & X^{l-1,q}  X^{l,1}  f_i(l)^{l-1,k}
\parallel_{L^2} \\ & \leq & \sum_{q+k=*-1} \parallel
X^{l-1,q}X^{l,1} X^{l-1,k} \cdots X^{r,i_{h_r}}f_i \parallel_{L^2}
\\&+& C \sum_{q+k=*-1}\Big(\sum_{y+z=k-1}
\parallel X^{l-1,q}X^{l,1}
X^{l-1,y}X^{l,1}V(l)^{l-1,z}
\parallel_{L^2}\Big)\\ &+& \sum_{q+k=*-1}\Big(\sum_{y+z=i_{(h_{s-1})-1}}\Big(
\sum_{s=l+1}^r
\parallel X^{l-1,q}X^{l,1} X^{l-1,k}X^{l,i_{h_l}} \cdots \\ &\cdots&
X^{s-1,y}X^{s,1}V(s)^{s-1,z}\parallel_{L^2}\Big)\Big)\,.
\end{eqnarray*}

Similarly, applying lemma \ref{L2} to the term $ X^{l-1,*}\cdots
X^{s-1,q}X^{s,1} f_i(s)^{s-1,k}$, we have:

\begin{eqnarray*}
\sum_{q+k=i_{(h_{s-1})-1}} & \parallel & X^{l-1,*}\cdots
X^{s-1,q}X^{s,1} f_i(s)^{s-1,k}\parallel_{L^2} \\& \leq &
\sum_{q+k=i_{(h_{s-1})-1}}\Big(\sum_{y+z=k-1} \Big( \sum_{s=l+1}^r
\parallel A_{ij}^{ \alpha
\beta}X^{l-1,*}\cdots  \\ &\cdots& X^{s-1,q}X^{s,1}X^{s-1,y}
X^{s,1}V(s)^{s-1,z}\parallel_{L^2}\Big)\Big) \\ &+&
\sum_{q+k=i_{(h_{s-1})-1}}\Big(\sum_{y+z=i_{(h_{p-1})-1}}\Big(\sum_{s=l+1}^r
\Big( \sum_{p=s+1}^r
\parallel A_{ij}^{ \alpha \beta}X^{l-1,*} \cdots \\ & \cdots &
X^{s-1,q}X^{s,1}X^{s-1,k}X^{s,i_{h_{s}}} \cdots
X^{p-1,y}X^{p,1}X^{p-1,z}X^{p,i_{h_{p}}}V(p+1)
\parallel_{L^2}\Big)\Big)\Big) \\&+& \sum_{q+k=i_{(h_{s-1})-1}} \Big(
\sum_{s=l+1}^r
\parallel X^{l-1,*} \cdots X^{s-1,q}X^{s,1}X^{s-1,k}X^{s,i_{h_{s}}}
\cdots X^{r,i_{h_r}}f_i \parallel_{L^2}\Big)\,.
\end{eqnarray*}

Substituting these estimates into the inequality derived in lemma
\ref{L1}, one arrives at the desired result.
\end{proof}

In order to use the above lemmas in the proofs of Theorems
$\ref{T1}$ and $\ref{T2}$, we need to modify them by grouping like
terms together.  Recall that differentiation along the last layer,
r, commutes with all other layers, so we can freely move these
derivatives around. However, since each of the other layers do not
commute we gain extra terms, called commutators, when we choose to
switch the order of differentiation. The following lemma will be
applied numerous times to these commutator terms that appear when
we group like terms together.

\begin{lemma}\label{L4}
Applying commutator properties, we have:
\begin{eqnarray}
X^{l-1,q}X^{l,1}X^{l-1,k}X^{l,i_{h_l}}V(l+1) &=&
X^{l-1,q+k}X^{l,i_{(h_l)+1}}V(l+1) \nonumber \\&+&
\sum_{s+t=q+k-1}X^{l-1,s}X^{2l-1,1}X^{l-1,t}X^{l,i_{h_l}}V(l+1)\,.
\end{eqnarray}
\end{lemma}

\begin{proof}
Begin with $X^{l-1,q}X^{l,1}X^{l-1,k}X^{l,i_{h_l}}V(l+1)$; we will
transfer $X^{l,1}$ to its like terms.  In order to do this, we
need to shift it "k" times to the right to get it past all k
derivatives in the "l-1" direction. To see how this process works,
first move $X^{l,1}$ just once to the right:
\begin{eqnarray*}
X^{l-1,q}X^{l,1}X^{l-1,k}X^{l,i_{h_l}}V(l+1) &=&
X^{l-1,q}X^{l-1,1}X^{l,1}X^{l-1,k-1}X^{l,i_{h_l}}V(l+1)\\&+&
X^{l-1,q}[X^{l,1},X^{l-1,1}]X^{l-1,k-1}X^{l,i_{h_l}}V(l+1)\\&=&
X^{l-1,q}X^{l-1,1}X^{l,1}X^{l-1,k-1}X^{l,i_{h_l}}V(l+1)\\&+&
X^{l-1,q}X^{2l-1,1}X^{l-1,k-1}X^{l,i_{h_l}}V(l+1)\,.
\end{eqnarray*}

From here it is clear to see that if we apply this same technique
"k-1" more times we eventually have the following:

\begin{eqnarray*}
X^{l-1,q}X^{l,1}X^{l-1,k}X^{l,i_{h_l}}V(l+1) &=&
X^{l-1,q}X^{l-1,k}X^{l,1}X^{l,i_{h_l}}V(l+1)
\\&+& X^{l-1,q}X^{2l-1,1}X^{l-1,k-1}X^{l,i_{h_l}}V(l+1) \\ &+&
X^{l-1,q}X^{l-1,1}X^{2l-1,1}X^{l-1,k-2}X^{l,i_{h_l}}V(l+1) +
\cdots
\\&+& X^{l-1,q}X^{l-1,k-1}X^{2l-1,1}X^{l,i_{h_l}}V(l+1)\\&=&
X^{l-1,q+k}X^{l,i_{(h_l)+1}}V(l+1)\\&+&
\sum_{s+t=q+k-1}X^{l-1,s}X^{2l-1,1}X^{l-1,t}X^{l,i_{h_l}}V(l+1)\,.
\end{eqnarray*}
\end{proof}

\subsection{Proof of Theorems 9 and 10}
\begin{proof}(Theorem \ref{T1})

Once again we let $*=i_{h_{l-1}}$.  We assume that $u^\beta \in
S^{1,2}_{loc}(\Omega)$ is a weak solution to the system
(\ref{1.1}), and we want to show that the following inequality is
finite:
\begin{eqnarray}\label{aa}
\nonumber \parallel X^{I_{h_{l-1}}}X^{I_{h_l}}& \cdots&
X^{I_{h_r}}u
\parallel_{S^{1,2}_{loc}}\\ &\leq &  \nonumber \Bigg[\parallel
X^{l-1,i_{h_{l-1}}}X^{I_{h_l}} \cdots X^{I_{h_r}}u^\beta
\parallel_{S^{1,2}_{loc}}+ \parallel
f^\al(l)^{l-1,*}\parallel_{L^2}\\ &+& \sum_{j=l-1}^r\Big(
\sum_{k=1}^{m_j} \parallel
X_{j,k}f^\al(l)^{l-1,*}\parallel_{L^2}\Big)+ \sum_{i=1}^m
\parallel f_i^\al(l)^{l-1,*}\parallel_{L^2}\\ &+& \nonumber
\sum_{i=1}^m \Big(\sum_{j=l-1}^r \Big(\sum_{k=1}^{m_j} \parallel
X_{j,k}f_i^\al(l)^{l-1,*}\parallel_{L^2}\Big)\Big) \Bigg]\,.
\end{eqnarray}

In order to show that $X^{I_{h_{l-1}}}X^{I_{h_l}} \cdots
X^{I_{h_r}}u$ is bounded from above, we need to show that each of
the terms on the right hand side of $\eqref{aa}$ is bounded from
above.  The terms that require the most work are
$X_{j,k}f^\al(l)^{l-1,*}$ and $ X_{j,k}f_i^\al(l)^{l-1,*}$, so we
will begin with these; this is where we use Lemmas \ref{L2} and
\ref{L3} since they essentially show us exactly what we are
looking at when we see $f^\al(l)^{l-1,*}$ and
$f_i^\al(l)^{l-1,*}$. Set

\begin{eqnarray*}
X_{j,k}f^\al(l)^{l-1,*} &=& X_{j,k}X^{l-1,*}X^{l,i_{h_l}} \cdots
X^{r,i_{h_r}}f^\al
+ X_{j,k}\Big(X^{l-1,*-1}X^{l,i_{(h_{l})+1}} \cdots X^{r,i_{h_r}}f_i^\al \\
&+& X_{j,k}\sum_{q+k=i_{(h_{s-1})-1}} \Big( \sum_{s=l+1}^r
X^{l-1,*} \cdots X^{s-1,q}X^{s,1}X^{s-1,k}X^{s,i_{h_{s}}} \cdots \\
&\cdots&
X^{r,i_{h_r}}f_i^\al \Big)\Big)\\
&+& P_1 + P_2 + P_3 + P_4 + P_5 + P_6 \,,
\end{eqnarray*}

where for simplicity we have let $P_j$, $j=1,\ldots,6$ equal the
following:

\begin{eqnarray*}
P_1&=&\sum_{q+k=*-1}  A_{ij}^{ \alpha \beta}
X_{j,k}X^{l-1,q}X^{l,1}X_jV(l)^{l-1,k}\,,\\
P_2&=&\sum_{s=l+1}^r \Big( \sum_{q+k=i_{(h_{s-1})-1}} A_{ij}^{
\alpha \beta}X_{j,k}X^{l-1,*}\cdots
X^{s-1,q}X^{s,1}X_jX^{s-1,k}X^{s,i_{h_s}} \cdots
X^{r,i_{h_r}}u^\beta\Big)\,,\\
P_3&=&C \sum_{q+k=*-1}\Big(\sum_{y+z=k-1}  X_{j,k}X^{l-1,q}X^{l,1}
X^{l-1,y}X^{l,1}V(l)^{l-1,z}\Big)\,,\\
P_4&=&\sum_{q+k=*-1}\Big(\sum_{y+z=i_{(h_{s-1})-1}}\Big(
\sum_{s=l+1}^r X_{j,k}X^{l-1,q}X^{l,1} X^{l-1,k}X^{l,i_{h_l}}
\cdots
X^{s-1,y}X^{s,1}V(s)^{s-1,z}\Big)\Big)\,,\\
P_5&=&\sum_{q+k=i_{(h_{s-1})-1}}\Big(\sum_{y+z=k-1} \Big(
\sum_{s=l+1}^r A_{ij}^{ \alpha \beta}X_{j,k}X^{l-1,*}\cdots
X^{s-1,q}X^{s,1}X^{s-1,y} X^{s,1}V(s)^{s-1,z}\Big)\Big)\,,\\
P_6&=&\sum_{q+k=i_{(h_{s-1})-1}}\Big(\sum_{y+z=i_{(h_{p-1})-1}}\Big(\sum_{s=l+1}^r
\Big( \sum_{p=s+1}^r
 A_{ij}^{ \alpha \beta}X_{j,k}X^{l-1,*} \cdots \\&\cdots&
X^{s-1,q}X^{s,1}X^{s-1,k}X^{s,i_{h_{s}}} \cdots
X^{p-1,y}X^{p,1}X^{p-1,z}X^{p,i_{h_{p}}}V(p+1)\Big)\Big)\Big)\,.
\end{eqnarray*}

Before we begin bounding each of the $P_j$ terms, recall from
$\eqref{aa}$ that for $X_{j,k}$ we are assuming $j \geq l-1$. If
$j> l-1$, then we always have fewer derivatives in the lowest
layer than when we began.  When this is the case, the estimate
$\eqref{T1eq1}$ in Theorem $\ref{T1}$ is satisfied. The only time
that we may not lessen the number of derivatives in the lowest
layer is if $j=l-1$, so this is what we will assume from here on.

When we apply Lemma \ref{L4} to the terms above that need
rearranging, things quickly get complicated.  However, the idea
behind this theorem is to not necessarily have to keep track of
each and every derivative, but instead to first count the total
derivatives on $u^\beta$ and second to count the number of
derivatives on $u^\beta$ in the lowest layer. We proceed by
writing an estimate for each of the $P_j$'s above in relation to
how many derivatives are attached to $u^\beta$ in each layer.  One
item worth noting is that when you "move" derivatives around using
lemma (\ref{L4}), you end up with numerous commutator terms, one
for each time you move a derivative that does not commute.  We do
not need to actually keep track of each of these terms; we just
note that when we commute two derivatives, we end up with one
derivative in a higher layer, thus lessening the total number of
derivatives. Therefore, we can group all of these terms together
in one collective term that we will call "commutator" and be
confident that this term has less total derivatives on $u$ than
what we started with.

$\underline{For P_1}$:
\begin{eqnarray*}
\parallel P_1 \parallel_{L^2} &\leq& \parallel \sum_{q+k=*-1}  A_{ij}^{
\alpha \beta} X_{j,k}X^{l-1,q}X^{l,1}X_jV(l)^{l-1,k}\parallel_{L^2} \\
&\leq& C\parallel
X_jX_{j,k}X^{l-1,*-1}X^{l,{h_l}+1}V(l+1)\parallel_{L^2}\\&+& C
\parallel \textrm{Commutator}\parallel_{L^{2}_{loc}}
\\ & \leq & C \parallel J_{h_{l-1}}J_{h_{l}} \cdots J_{h_r} u^\beta
\parallel_{S^{1,2}_{loc}}
\end{eqnarray*}
where either $\sum_{i=l-1}^r J_{h_i} < \sum_{i=l-1}^r I_{h_i}$ (as
is the case in the commutator term) or $\sum_{i=l-1}^r J_{h_i} =
\sum_{i=l-1}^r I_{h_i}$ so that $|J_{h_{l-1}}| < |I_{h_{l-1}}|$
and $|J_{h_{l}}|
> |I_{h_{l}}|$ with $J_{h_k} = I_{h_k}$ elsewhere.  Thus, $P_1$ satisfies
estimate $\eqref{T1eq1}$ in the statement of the theorem.

$\underline{For P_2}$:
\begin{eqnarray*}
\parallel P_2 \parallel_{L^2} & \leq & \sum_{s=l+1}^r \Big(
\sum_{q+k=i_{(h_{s-1})-1}} \parallel A_{ij}^{ \alpha
\beta}X_{j,k}X^{l-1,*}\cdots
X^{s-1,q}X^{s,1}X_jX^{s-1,k}X^{s,i_{h_s}}\\ &\cdots&
X^{r,i_{h_r}}u^\beta \parallel_{L^2} \Big)\\ &\leq& C
\sum_{s=l+1}^r \parallel X_jX_{j,k}X^{l-1,*}\cdots
X^{s-1,{h_{s-1}-1}}X^{{s,{h_s}+1}} \cdots  X^{r,i_{h_r}}u^\beta
\parallel_{L^2} \\ &+& \parallel \textrm{Commutator}\parallel_{L^{2}_{loc}}
\\
& \leq & C \parallel J_{h_{l-1}}J_{h_{l}} \cdots J_{h_r} u^\beta
\parallel_{S^{1,2}_{loc}}
\end{eqnarray*}
where  either $\sum_{i=l-1}^r J_{h_i} < \sum_{i=l-1}^r I_{h_i}$
(as is the case in the commutator term) or $\sum_{i=l-1}^r J_{h_i}
= \sum_{i=l-1}^r I_{h_i}$ such that one of three things happens:
Either $J_{h_{l-1}} = I_{h_{l-1}}$ and there exists $\beta_i
> h_{l-1}$, $i \geq 1$, such that $J_{\beta_i} > I_{\beta_i}$ with
$J_{h_k} \leq I_{h_k}$ for every $h_k \neq h_{l-1},\beta_i$ or
$J_{h_{l-1}} = I_{h_{l-1}}$ with $J_{h_{s-1}} < I_{h_{s-1}}$
(actually $= I_{h_{s-1}} -1$), $J_{h_{s}} > I_{h_{s}}$ (actually
$= I_{h_{s}} +1$) and $J_{h_k} \leq I_{h_k}$ for every $k \neq
l-1, s-1,s$ or $J_{h_{l-1}} = I_{h_{l-1}}$ and there exists
$\beta_i
> h_{l-1}$, $i \geq 1$, such that $J_{\beta_i} > I_{\beta_i}$ with
$J_{h_k} \leq I_{h_k}$ for every $h_k \neq h_{l-1},\beta_i$.
Thus, $P_2$ satisfies estimate $\eqref{T1eq1}$ in the statement of
the theorem.

$\underline{For P_3}$:
\begin{eqnarray*}
\parallel P_3 \parallel_{L^2} & \leq & C
\sum_{q+k=*-1}\Bigg(\sum_{y+z=k-1} \parallel
X_{j,k}X^{l-1,q}X^{l,1}
X^{l-1,y}X^{l,1}V(l)^{l-1,z}\parallel_{L^2}\Bigg) \\
& \leq & C \parallel J_{h_{l-1}}J_{h_{l}} \cdots J_{h_r} u^\beta
\parallel_{S^{1,2}_{loc}}
\end{eqnarray*}

where either $\sum_{i=l-1}^r J_{h_i} = \sum_{i=l-1}^r I_{h_i}$ so
that $J_{h_{l-1}} < I_{h_{l-1}}$ (actually $= I_{h_{l-1}} -2$),
$J_{h_l}
> I_{h_l}$ (actually $ = I_{h_l} +2 $), and $J_{h_k} = I_{h_k}$ for every $k
\neq l-1,l$ or $\sum_{i=l-1}^r J_{h_i} < \sum_{i=l-1}^r I_{h_i}$
with $J_{h_{l-1}}$ always at least 2 less than $I_{h_{l-1}}$ and
there exists $\beta_i > h_{l-1}$ for $i \geq 1$ such that
$J_{\beta_i} > I_{\beta_i}$ with $J_{h_k} \leq I_{h_k}$ for all
$h_k \neq h_{l-1}, \beta_i$.  Thus, $P_3$ satisfies estimate
$\eqref{T1eq1}$ in the statement of the theorem.

$\underline{For P_4}$:
\begin{eqnarray*}
\parallel P_4 \parallel_{L^2} & \leq &
\sum_{q+k=*-1}\Bigg(\sum_{y+z=i_{(h_{s-1})-1}}\\ &\Big(&
\sum_{s=l+1}^r
\parallel X_{j,k}X^{l-1,q}X^{l,1} X^{l-1,k}X^{l,i_{h_l}} \cdots
X^{s-1,y}X^{s,1}V(s)^{s-1,z} \parallel_{L^2} \Big)\Bigg) \\
& \leq & C \parallel J_{h_{l-1}}J_{h_{l}} \cdots J_{h_r} u^\beta
\parallel_{S^{1,2}_{loc}}
\end{eqnarray*}
where either $\sum_{i=l-1}^r J_{h_i} = \sum_{i=l-1}^r I_{h_i}$ so
that $J_{h_{l-1}} < I_{h_{l-1}}$ (actually $= I_{h_{l-1}} -1$),
$J_{h_l}
> I_{h_l}$ (actually $ = I_{h_l} +1 $), $J_{h_{s-1}} < I_{h_{s-1}}$ (actually $=
I_{h_{s-1}} -1$), $J_{h_s} > I_{h_s}$ (actually $ = I_{h_s} +1 $),
and $J_{h_k} = I_{h_k}$ for every $k \neq l-1,l,s-1,s$ or
$\sum_{i=l-1}^r J_{h_i} < \sum_{i=l-1}^r I_{h_i}$ with
$J_{h_{l-1}}$ and $J_{h_{s-1}}$ always at least 1 less than
$I_{h_{l-1}}$ and $I_{h_{s-1}}$, respectively, and then there
exists $\beta_i > h_{l-1}$ and/or $ \beta_i > h_{s-1}$ for $i \geq
1$ such that $J_{\beta_i}
> I_{\beta_i}$ with $J_{h_k} \leq I_{h_k}$ for all $h_k \neq h_{l-1}, h_{s-1},
\beta_i$.  Thus, $P_4$ satisfies estimate $\eqref{T1eq1}$ in the
statement of the theorem.

$\underline{For P_5}$:
\begin{eqnarray*}
\parallel P_5 \parallel_{L^2} & \leq & \sum_{q+k+z=i_{(h_{s-1})-2}}
\Bigg( \sum_{s=l+1}^r \parallel A_{ij}^{ \alpha
\beta}X_{j,k}X^{l-1,*}\cdots X^{s-1,q}X^{s,1}X^{s-1,y}
X^{s,1}V(s)^{s-1,z} \parallel_{L^2}\Bigg)\\
& \leq & C \parallel J_{h_{l-1}}J_{h_{l}} \cdots J_{h_r} u^\beta
\parallel_{S^{1,2}_{loc}}
\end{eqnarray*}
where either $\sum_{i=l-1}^r J_{h_i} = \sum_{i=l-1}^r I_{h_i}$ so
that $J_{h_{l-1}} = I_{h_{l-1}}$, $J_{h_{s-1}} < I_{h_{s-1}}$
(actually $= I_{h_{s-1}} -2$), $J_{h_s} > I_{h_s}$ (actually $ =
I_{h_s} +2 $), and $J_{h_k} = I_{h_k}$ for every $k \neq s-1,s$ or
$\sum_{i=l-1}^r J_{h_i} < \sum_{i=l-1}^r I_{h_i}$ with
$J_{h_{s-1}}$ always at least 2 less than $I_{h_{s-1}}$, and then
there exists $ \beta_i > h_{s-1}$ for $i \geq 1$ such that
$J_{\beta_i}
> I_{\beta_i}$ with $J_{h_k} \leq I_{h_k}$ for all $h_k \neq h_{s-1},
\beta_i$. Thus, $P_5$ satisfies estimate $\eqref{T1eq1}$ in the
statement of the theorem.

$\underline{For P_6}$:
\begin{eqnarray*}
\parallel P_6 \parallel_{L^2} & \leq &
\sum_{q+k=i_{(h_{s-1})-1}}\Bigg(\sum_{y+z=i_{(h_{p-1})-1}}\Bigg(\sum_{s=l+1}^r
\Bigg( \sum_{p=s+1}^r \parallel A_{ij}^{ \alpha
\beta}X_{j,k}X^{l-1,*}\\ &\cdots&
X^{s-1,q}X^{s,1}X^{s-1,k}X^{s,i_{h_{s})}} \cdots \\&\cdots&
X^{p-1,y}X^{p,1}X^{p-1,z}X^{p,i_{h_{p}}}V(p+1) \parallel_{L^2}\Bigg)\Bigg)\Bigg)\\
& \leq & C \parallel J_{h_{l-1}}J_{h_{l}} \cdots J_{h_r} u^\beta
\parallel_{S^{1,2}_{loc}}
\end{eqnarray*}
where either $\sum_{i=l-1}^r J_{h_i} = \sum_{i=l-1}^r I_{h_i}$ so
that $J_{h_{s-1}} < I_{h_{s-1}}$ (actually $= I_{h_{s-1}} -1$),
$J_{h_s}
> I_{h_s}$ (actually $ = I_{h_s} +1 $), $J_{h_{p-1}} < I_{h_{p-1}}$ (actually $=
I_{h_{p-1}} -1$), $J_{h_p} > I_{h_p}$ (actually $ = I_{h_p} +1 $),
and $J_{h_k} = I_{h_k}$ for every $k \neq s-1,s,p-1,p$ or
$\sum_{i=l-1}^r J_{h_i} < \sum_{i=l-1}^r I_{h_i}$ with
$J_{h_{s-1}}$ and $J_{h_{p-1}}$ always at least 1 less than
$I_{h_{s-1}}$ and $I_{h_{p-1}}$, respectively, and then there
exists $\beta_i > h_{s-1}$ and/or $ \beta_i > h_{p-1}$ for $i \geq
1$ such that $J_{\beta_i}
> I_{\beta_i}$ with $J_{h_k} \leq I_{h_k}$ for all $h_k \neq h_{s-1}, h_{p-1},
\beta_i$.  Thus, $P_6$ satisfies estimate $\eqref{T1eq1}$ in the
statement of the theorem.

Combining these estimates we have the following:

\begin{eqnarray}\label{r}
\nonumber \parallel X_{j,k}f^\al(l)^{l-1,*}\parallel_{L^2} & \leq
&
\parallel X_{j,k}X^{l-1,*}X^{l,i_{h_l}} \cdots X^{r,i_{h_r}}f^\al
\parallel_{L^2}\\ &+& \nonumber
\parallel X_{j,k}X^{l-1,*-1}X^{l,i_{(h_l)+1}} \cdots
X^{r,i_{h_r}}f_i^\al \parallel_{L^2} \\
&+&  \nonumber \parallel X_{j,k}\sum_{q+k=i_{(h_{s-1})-1}}
\Big(\sum_{s=l+1}^r X^{l-1,*} \cdots \\ & \cdots & \nonumber
X^{s-1,q}X^{s,1}X^{s-1,k}X^{s,i_{h_{s}}} \cdots
X^{r,i_{h_r}}f_i^\al \parallel_{L^2} \Big)\\
&+& \parallel P_1 \parallel_{L^2} + \parallel P_2\parallel_{L^2} +
\parallel P_3\parallel_{L^2} + \parallel P_4\parallel_{L^2} +
\parallel P_5\parallel_{L^2}+ \parallel P_6\parallel_{L^2} \\
& \leq & \nonumber \parallel X_{j,k}X^{l-1,*}X^{l,i_{h_l}} \cdots
X^{r,i_{h_r}}f^\al \parallel_{L^2}\\& +& \nonumber
\parallel X_{j,k}X^{l-1,*-1}X^{l,i_{(h_l)+1}} \cdots
X^{r,i_{h_r}}f_i^\al \parallel_{L^2} \\
&+& \nonumber \parallel X_{j,k}\sum_{q+k=i_{(h_{s-1})-1}} \Big(
\sum_{s=l+1}^r X^{l-1,*} \cdots\\ &\cdots& \nonumber
X^{s-1,q}X^{s,1}X^{s-1,k}X^{s,i_{h_{s}}} \cdots
X^{r,i_{h_r}}f_i^\al
\parallel_{L^2}\Big)\\ &+& \nonumber \parallel J_{h_{l-1}}J_{h_l} \cdots
J_{h_r} u^\beta
\parallel_{S^{1,2}_{loc}}\,.
\end{eqnarray}

Next, set
\begin{eqnarray*}
X_{j,k}f_i^\al(l)^{l-1,*}&=& X_{j,k}X^{l-1,*}X^{l,i_{h_l}} \cdots
X^{r,i_{h_r}}f_i^\al
\\&+& \sum_{s=l+1}^r \Big( \sum_{q+k=i_{(h_{s-1})-1}}A_{ij}^{ \alpha
\beta}X_{j,k}X^{l-1,*} \cdots X^{s-1,q}X^{s,1}V(s)^{s-1,k} \Big)\\
&+& \sum_{q+k=*-1} A_{ij}^{ \alpha \beta}
X_{j,k}X^{l-1,q}X^{l,1}V(l)^{l-1,k}
\\ &=& X_{j,k}X^{l-1,*}X^{l,i_{h_l}}
\cdots X^{r,i_{h_r}}f_i^\al \\ &+& Q_1 + Q_2\,,
\end{eqnarray*}

where we let
\begin{eqnarray*}
Q_1 &=&\sum_{s=l+1}^r\Big( \sum_{q+k=i_{(h_{s-1})-1}}A_{ij}^{
\alpha
\beta}X_{j,k}X^{l-1,*} \cdots X^{s-1,q}X^{s,1}V(s)^{s-1,k}\Big)\,,\\
Q_2&=&\sum_{q+k=*-1} A_{ij}^{ \alpha \beta}
X_{j,k}X^{l-1,q}X^{l,1}V(l)^{l-1,k}\,.
\end{eqnarray*}

Proceeding as we did before by applying the results of Lemma
\ref{L4}, we have the following:

$\underline{For Q_1}$:

\begin{eqnarray*}
\parallel Q_1 \parallel_{L^2} & \leq & \sum_{s=l+1}^r
\Bigg(\sum_{q+k=i_{(h_{s-1})-1}}\parallel A_{ij}^{ \alpha
\beta}X_{j,k}X^{l-1,*} \cdots X^{s-1,q}X^{s,1}V(s)^{s-1,k}
\parallel_{L^2}\Bigg) \\ & \leq & C \parallel J_{h_{l-1}}J_{h_{l}} \cdots
J_{h_r} u^\beta
\parallel_{S^{1,2}_{loc}}
\end{eqnarray*}
where either $\sum_{i=l-1}^r J_{h_i} = \sum_{i=l-1}^r I_{h_i}$
such that $J_{h_{l-1}}=I_{h_{l-1}}$, $J_{h_{s-1}} < I_{h_{s-1}}$
(actually $= I_{h_{s-1}} -1$), $J_{h_s}
> I_{h_s}$ (actually $ = I_{h_s} +1 $), and $J_{h_k} = I_{h_k}$ for every $k
\neq l-1, s-1, s$ or $\sum_{i=l-1}^r J_{h_i} < \sum_{i=l-1}^r
I_{h_i}$ with $J_{h_{l-1}}=I_{h_{l-1}}$, $J_{h_{s-1}}$ always at
least 1 less than $I_{h_{s-1}}$, and then there exists $\beta_i >
h_{s-1}$ for $i \geq 1$ such that $J_{\beta_i}
> I_{\beta_i}$ with $J_{h_k} \leq I_{h_k}$ for all $h_k \neq h_{l-1}, h_{s-1},
\beta_i$.  Thus, $Q_1$ satisfies estimate $\eqref{T1eq1}$ in the
statement of the theorem.

$\underline{For Q_2}$:

\begin{eqnarray*}
\parallel Q_2 \parallel_{L^2} & \leq & \sum_{q+k=*-1} \parallel
A_{ij}^{ \alpha \beta}
X_{j,k}X^{l-1,q}X^{l,1}V(l)^{l-1,k}\parallel_{L^2}\\ & \leq & C
\parallel J_{h_{l-1}}J_{h_{l}} \cdots J_{h_r} u^\beta
\parallel_{S^{1,2}_{loc}}
\end{eqnarray*}
where either $\sum_{i=l-1}^r J_{h_i} = \sum_{i=l-1}^r I_{h_i}$
such that $J_{h_{l-1}} < I_{h_{l-1}}$ (actually $= I_{h_{l-1}}
-1$), $J_{h_l}
> I_{h_l}$ (actually $ = I_{h_l} +1 $), and $J_{h_k} = I_{h_k}$ for every $k
\neq l-1, l$ or $\sum_{i=l-1}^r J_{h_i} < \sum_{i=l-1}^r I_{h_i}$
with $J_{h_{l-1}} < I_{h_{l-1}}$, and then there exists $\beta_i
> h_{l-1}$ for $i \geq 1$ such that $J_{\beta_i}
> I_{\beta_i}$ with $J_{h_k} \leq I_{h_k}$ for all $h_k \neq h_{l-1},
\beta_i$. Thus, $Q_2$ satisfies estimate $\eqref{T1eq1}$ in the
statement of the theorem.

Combining these estimates, we have the following:

\begin{eqnarray}\label{s}
\nonumber \parallel X_{j,k}f_i^\al(l)^{l-1,*}\parallel_{L^2}& \leq
&
\parallel X_{j,k}X^{l-1,*}X^{l,i_{h_l}} \cdots
X^{r,i_{h_r}}f_i^\al
\parallel_{L^2}
\\&+& \nonumber \sum_{s=l+1}^r\Big( \sum_{q+k=i_{(h_{s-1})-1}} C \parallel
X_{j,k}X^{l-1,*} \cdots X^{s-1,q}X^{s,1}V(s)^{s-1,k}\parallel_{L^2}\Big)\\
&+& \nonumber \sum_{q+k=*-1} \parallel
X_{j,k}X^{l-1,q}X^{l,1}V(l)^{l-1,k}\parallel_{L^2}
\\ & \leq & \parallel X_{j,k}X^{l-1,*}X^{l,i_{h_l}}
\cdots X^{r,i_{h_r}}f_i^\al\parallel_{L^2} \\ &+& \nonumber
\parallel Q_1
\parallel_{L^2} + \parallel Q_2 \parallel_{L^2}
\\ & \leq & \nonumber \parallel X_{j,k}X^{l-1,*}X^{l,i_{h_l}}
\cdots X^{r,i_{h_r}}f_i^\al\parallel_{L^2} \\ &+&  \nonumber C
\parallel J_{h_{l-1}}J_{h_{l}} \cdots J_{h_r} u^\beta
\parallel_{S^{1,2}_{loc}}\,.
\end{eqnarray}

What we have done is bound the terms $X_{j,k}f^\al(l)^{l-1,*}$ and
$X_{j,k}f_i^\al(l)^{l-1,*}$ from above by counting and keeping
track of the derivatives in each step. We will also need bounds
for $f^\al(l)^{l-1,*}$ and $f_i^\al(l)^{l-1,*}$, but since we are
counting derivatives, all of our estimates are just one less than
what we calculated above.  It is clear, then, that the total
number of derivatives in these terms is less than what we started
with so they satisfy the estimate $\eqref{T1eq1}$ in Theorem
$\ref{T1}$.

If we look at the terms that involve differentiation on $f^\al$
and $f_i^\al$, we can group these together and collectively name
them $\tilde{f^\al}$ and $\tilde{f_i^\al}$. Since our original
$f^\al$ and $f_i^\al$ are in $C^\infty$, we know that these are
bounded:

\begin{eqnarray}\label{t}
 \parallel \tilde{f}^\alpha
\parallel_{L^2} &\leq & \parallel X^{l-1,*}X^{l,i_{h_l}} \cdots
X^{r,i_{h_r}}f^\alpha
\parallel_{L^2} \notag \\ &+&  \sum_{j=l-1}^r \sum_{k=1}^{m_j}\parallel
X_{j,k}X^{l-1,*}X^{l,i_{h_l}} \cdots X^{r,i_{h_r}}f^\alpha
\parallel_{L^2},\notag
\\ \textrm{and}\notag
\\ \parallel \tilde{f_i}^\alpha \parallel_{L^2}
& \leq&
\parallel X^{l-1,*}X^{l,i_{h_l}} \cdots X^{r,i_{h_r}}f_i^\alpha
\parallel_{L^2} + \parallel X^{l-1,*-1}X^{l,i_{(h_{l})+1}} \cdots
X^{r,i_{h_r}}f_i^\alpha
\parallel_{L^2}\notag \\ &+&  \sum_{q+k=i_{(h_{s-1})-1}} \Big(\sum_{s=l+1}^r
\parallel X^{l-1,*} \cdots X^{s-1,q}X^{s,1}X^{s-1,k}X^{s,i_{h_{s}}}
\cdots X^{r,i_{h_r}}f_i^\alpha \parallel_{L^2}\Big)\notag \\
&+& \sum_{j=l-1}^r \sum_{k=1}^{m_j} X_{j,k} \Big(\parallel
X^{l-1,*}X^{l,i_{h_l}} \cdots X^{r,i_{h_r}}f_i^\alpha
\parallel_{L^2}\notag \\ &&+  \parallel X^{l-1,*-1}X^{l,i_{(h_l)+1}}
\cdots X^{r,i_{h_r}}f_i^\alpha
\parallel_{L^2} \\ &+&  \sum_{q+k=i_{(h_{s-1})-1}} \Big( \sum_{s=l+1}^r
\parallel X^{l-1,*} \cdots X^{s-1,q}X^{s,1}X^{s-1,k}X^{s,i_{h_{s}}}
\cdots X^{r,i_{h_r}}f_i^\alpha
\parallel_{L^2}\Big)\Big)\notag \,.
\end{eqnarray}

The following estimate follows from $\eqref{r}$, $\eqref{s}$, and
$\eqref{t}$:

\begin{eqnarray*}
\parallel X^{I_{h_{l-1}}}X^{I_{h_l}}& \cdots&
X^{I_{h_r}}u^\beta \parallel_{S^{1,2}_{loc}}\\ &\leq &
\sum_{i=1}^m
\parallel \tilde{f_i}^\alpha \parallel_{L^2} +\parallel
\tilde{f}^\alpha
\parallel_{L^2}\\ &+& C
\parallel J_{h_{l-1}}J_{h_{l}} \cdots J_{h_r} u^\beta
\parallel_{S^{1,2}_{loc}}\,,
\end{eqnarray*}
such that either $$ \sum_{i=l-1}^r J_{h_i} < \sum_{i=l-1}^r
I_{h_i}$$ or
$$ \sum_{i=l-1}^r J_{h_i} = \sum_{i=l-1}^r I_{h_i}\,,$$ so that there exists
$\beta_i > h_{l-1}$ such that $J_{\beta_i} > I_{\beta_i}$ and
$J_{h_k} \leq I_{h_k}$ for every $h_k \neq h_{l-1}, \beta_i$.
\end{proof}

\begin{proof}(Theorem \ref{T2})
We want an estimate for $\parallel
X^{I_{h_{l-1}}}X^{I_{h_r}}u^\beta
\parallel_{S^{1,2}_{loc}}$, so first consider what we have if we
differentiate (\ref{1.1}) in layer $r$. Since this layer commutes
with everything, applying $X^{I_{h_r}}$ simply gives us

\begin{eqnarray}\label{21}
\sum_{ \beta = 1}^N \sum _{i,j=1}^m X_i (A_{ij}^{ \alpha \beta}
X_j (X^{I_r}u^ {\beta}) + X^{I_{h_r}}f_i^{ \alpha}) =
X^{I_{h_r}}f^{ \alpha}\quad.
\end{eqnarray}

It is easy to see that $X^{I_{h_r}}u^\beta$ is indeed in
${S^{1,2}_{loc}}$, so next we need to show that we can
differentiate (\ref{21}) in the direction $l-1$:

By definition, we have that $$f_i^\alpha(r)^{l-1,*} \, = \,
A_{ij}^{ \alpha \beta}[X_{l-1,*},X_j]V(r)^{l-1,*-1}\, +\,
X_{l-1,*}f_i^\alpha(r)^{l-1,*-1}\,.$$  We can make use of the
proof of lemma (\ref{L2}) and our commutator result in lemma
(\ref{L4}) to see that
\begin{eqnarray*}
f_i^\alpha(r)^{l-1,*} &=& \sum_{q+k=*-1}A_{ij}^{ \alpha
\beta}X^{l-1,q}X^{l,1}V(r)^{l-1,k} + X^{l-1,*}f_i^
\alpha(r)^{l-1,0}\\&=& \sum_{q+k=*-1}A_{ij}^{ \alpha
\beta}X^{l-1,q}X^{l,1}X^{l-1,k}X^{I_{h_r}}u^ \beta +
X^{l-1,*}X^{I_{h_r}}f_i^ \alpha \\ &=& A_{ij}^{ \alpha
\beta}X^{l-1,*-1}X^{l,1}X^{I_{h_r}}u^ \beta \\ &+&
\sum_{s+t=*-2}A_{ij}^{\alpha
\beta}X^{l-1,s}X^{2l-1,1}X^{l-1,t}X^{I_{h_r}}u^ \beta \\ &+&
X^{l-1,*}X^{I_{h_r}}f_i^ \alpha\,.
\end{eqnarray*}

Continuing to apply Lemma \ref{L4} indefinitely in order to group
all like terms together, one obtains the following:

\begin{eqnarray*}
\parallel f_i^\alpha(r)^{l-1,*} \parallel_{L^2} & \leq &
C \parallel X^{J_{l-1}}X^{J_l} \cdots X^{J_r}u^\beta
\parallel_{L^2} +
\parallel X^{l-1,*}X^{I_{h_r}}f_i^ \alpha \parallel_{L^2}\,.
\end{eqnarray*}

Counting derivatives at this step, we see that $\sum_{i=l-1}^r
J_{h_i} < \sum_{i=l-1}^r I_{h_i}$ with $J_{h_{l-1}}$ always being
at least 2 less than $I_{h_{l-1}}$, $J_{h_k} > I_{h_k}$ for every
other $k$. The key here is that the derivatives are definitely at
least 2 less in the lowest layer; even though we are adding some
derivatives to the right of this layer, we have that the sum of
the derivatives in all layers is always less than what we started
with.

Next, we need to apply this same technique to
$f^\alpha(r)^{l-1,*}$.  Referring back to the definition and using
lemma (\ref{L1}) we have:
\begin{eqnarray*}
f^\alpha(r)^{l-1,*} \,&=& \, X_{l-1,*}f^\alpha(r)^{l-1,*-1}
\\ &+& \, [X_i,X_{l-1,*}] \big(A_{ij}^{ \alpha \beta} X_j
V(r)^{l-1,*-1} \, + \,f_i^\alpha(r)^{l-1,*-1}\big)\\
&=& X^{l-1,*}f^\alpha(r)^{l-1,0} +
\sum_{q+k=*-1}A_{ij}^{ \alpha \beta} X^{l-1,q}X^{l,1}X_jV(r)^{l-1,k}\\
&+& \sum_{q+k=*-1}X^{l-1,q}X^{l,1}f_i^\alpha(r)^{l-1,k}\\&=&
X^{l-1,*}X^{I_{h_r}}f^\alpha + \sum_{q+k=*-1}A_{ij}^{ \alpha
\beta} X^{l-1,q}X^{l,1}X_jX^{l-1,k}X^{I_{h_r}}u^ \beta \\ &+&
\sum_{q+k=*-1}X^{l-1,q}X^{l,1}f_i^\alpha(r)^{l-1,k}\,.
\end{eqnarray*}

From here the idea is the same as before.  Apply Lemma \ref{L4} to
$$\sum_{q+k=*-1}A_{ij}^{ \alpha \beta}
X^{l-1,q}X^{l,1}X_jX^{l-1,k}X^{I_{h_r}}u^ \beta$$ indefinitely in
order to group all like terms together.  However, what is
important and can actually be seen already is that we will always
have at least 2 derivatives less in the $l-1$ layer.  We may once
again add derivatives to the right of this layer, but we will
always have that the sum of all derivatives is less than what we
started with.

We still need to check
$\sum_{q+k=*-1}X^{l-1,q}X^{l,1}f_i^\alpha(r)^{l-1,k}$ to ensure it
follows this same mode of thought.  Applying previous calculations
we see that
\begin{eqnarray*}
\sum_{q+k=*-1}X^{l-1,q}X^{l,1}f_i^\alpha(r)^{l-1,k} &=& A_{ij}^{
\alpha
\beta}\sum_{q+k=*-1}X^{l-1,q}X^{l,1}X^{l-1,k-1}X^{l,1}X^{I_{h_r}}u^
\beta
\\ &+& A_{ij}^{\alpha
\beta}\sum_{q+k=*-1}X^{l-1,q}X^{l,1}\Big(\sum_{s+t=k-2}X^{l-1,s}X^{2l-1,1}X^{l-1,t}X^{I_{h_r}}u^
\beta\Big) \\ &+&
\sum_{q+k=*-1}X^{l-1,q}X^{l,1}X^{l-1,k}X^{I_{h_r}}f_i^ \alpha\,.
\end{eqnarray*}
Upon closer inspection of this equation, it is clear that we can
argue the same as before. Counting the derivatives in every layer
we see that the sum is at least 2 less than what we started with.
We also have at least 3 less in the lower layer, but we are again
adding derivatives to the higher layers.  The key once again is
that we are keeping less derivatives than what we started with and
that we have less in the lower layer.

Hence, we can say that the following inequality holds:

\begin{eqnarray*}
\parallel f^\alpha(r)^{l-1,*} \parallel_{L^2} & \leq &
C \Big[\parallel X^{J_{h_{l-1}}}X^{J_{h_l}} \cdots
X^{J_{h_r}}\parallel_{L^2} + \parallel
X_jX^{J_{h_{l-1}}}X^{J_{h_l}} \cdots X^{J_{h_r}}\parallel_{L^2}
\\ &+&
\parallel X^{l-1,*}X^{I_{h_r}}f^ \alpha \parallel_{L^2} +
\sum_{q+k=*-1}\parallel X^{l-1,q}X^{l,1}X^{l-1,k}X^{I_{h_r}}f_i^
\alpha
\parallel_{L^2}\Big]\,,
\end{eqnarray*}
where we have that the $\sum_{i=l-1}^r J_{h_i} < \sum_{i=l-1}^r
I_{h_i}$ with $J_{h_{l-1}}$ always being at least 2 less than
$I_{h_{l-1}}$ and $J_{h_k} > I_{h_k}$ for every other $k$.

We have now done the bulk of the calculations that we need in
order to gain an estimate on $X^{I_{h_{l-1}}}X^{I_{h_r}}u^\beta$.
Putting it all together we have the following:

\begin{eqnarray*}
\parallel X^{I_{h_{l-1}}}X^{I_{h_r}}u^\beta \parallel_{S^{1,2}_{loc}}&
\leq & C \bigg( \parallel X^{l-1,*}X^{I_{h_r}}u^\beta
\parallel_{S^{1,2}_{loc}}
+ \parallel f^\alpha(r)^{l-1,*}\parallel_{L^2}\\ &+&
\sum_{j=l-1}^r \sum_{k=1}^{m_j} \parallel
X_{j,k}f^\alpha(r)^{l-1,*}\parallel_{L^2}+ \sum_{i=1}^m
\parallel f_i^\alpha(r)^{l-1,*}\parallel_{L^2}\\ &+&
\sum_{i=1}^m \sum_{j=l-1}^r \sum_{k=1}^{m_j} \parallel
X_{j,k}f_i^\alpha(r)^{l-1,*}\parallel_{L^2} \bigg)\\ & \leq& C
\bigg(\parallel X^{J_{h_{l-1}}}X^{J_{h_l}} \cdots
X^{J_{h_r}}u^\beta
\parallel_{S^{1,2}_{loc}} + \sum_{i=1}^m \parallel
X^{l-1,*}X^{I_{h_r}}f_i^ \alpha \parallel_{L^2} \\ &+&
\parallel X^{l-1,*}X^{I_{h_r}}f^ \alpha \parallel_{L^2} +
\sum_{i=1}^m \Big(\sum_{q+k=*-1}\parallel
X^{l-1,q}X^{l,1}X^{l-1,k}X^{I_{h_r}}f_i^ \alpha
\parallel_{L^2}\Big) \bigg)\\ &+&  C\sum_{j=l-1}^r \sum_{k=1}^{m_j}
\bigg(\parallel X_{j,k}X^{J_{h_{l-1}}}X^{J_{h_l}} \cdots
X^{J_{h_r}}u^\beta
\parallel_{S^{1,2}_{loc}}\\ & +& \sum_{i=1}^m \parallel
X_{j,k}X^{l-1,*}X^{I_{h_r}}f_i^ \alpha \parallel_{L^2} \\ &+&
\parallel X_{j,k}X^{l-1,*}X^{I_{h_r}}f^ \alpha \parallel_{L^2}
\\&+&
\sum_{i=1}^m \Big(\sum_{q+k=*-1}\parallel
X_{j,k}X^{l-1,q}X^{l,1}X^{l-1,k}X^{I_{h_r}}f_i^ \alpha
\parallel_{L^2}\Big) \bigg)
\\ & \leq & C \bigg(\parallel X^{J_{h_{l-1}}}X^{J_{h_l}} \cdots
X^{J_{h_r}}u^\beta
\parallel_{S^{1,2}_{loc}} + \sum_{i=1}^m \parallel
\tilde{f_i}^ \alpha \parallel_{L^2} + \parallel \tilde{f}^ \alpha
\parallel_{L^2}\bigg)
\end{eqnarray*}

such that $\sum_{i=l-1}^r J_{h_i} < \sum_{i=l-1}^r I_{h_i}$ with
$J_{h_{l-1}}$ always being at least 1 less than $I_{h_{l-1}}$ and
$J_{h_k} \geq I_{h_k}$ for every other $k$,

where we set
\begin{eqnarray*}
\tilde{f} =X^{l-1,*}X^{I_{h_r}}f^ \alpha
+X_{j,k}X^{l-1,*}X^{I_{h_r}}f^ \alpha
\end{eqnarray*}

and

\begin{eqnarray*}
\tilde{f_i} &=&X^{l-1,*}X^{I_{h_r}}f_i^ \alpha +\sum_{q+k=*-1}
X^{l-1,q}X^{l,1}X^{l-1,k}X^{I_{h_r}}f_i^ \alpha
\\&+&X_{j,k}X^{l-1,*}X^{I_{h_r}}f_i^ \alpha +
\sum_{q+k=*-1}X_{j,k}X^{l-1,q}X^{l,1}X^{l-1,k}X^{I_{h_r}}f_i^
\alpha\,.
\end{eqnarray*}

\end{proof}

Iterating Theorems \ref{T1} and \ref{T2} we have the following
result:

\begin{cor}\label{c1}
Let $u \in S^{1,2}_{loc}(\Omega)$ be a weak solution to
(\ref{1.1}) such that the hypothesis of Theorems \ref{T1} and
\ref{T2} hold. Then the following is true:
\begin{eqnarray*}
\parallel X^{I_1} \cdots X^{I_r}
u\parallel_{S^{1,2}_{loc}(B(0,1))} &\leq& C \Big( \parallel
X^{I_r}u \parallel_{S^{1,2}_{loc}(B(0,2))} + \parallel \tilde f
\parallel_{L^2(B(0,2))} + \parallel \tilde {f_i}
\parallel_{L^2(B(0,2))} \Big)\,,
\end{eqnarray*}
where $\tilde{f}$, $\tilde{f_i}$ represent combinations of
derivatives on $f$, $f_i$, respectively.
\end{cor}

\begin{proof}(Main Theorem)

Iterating Theorem \ref{A}, we have the following estimate
\begin{eqnarray*}
\parallel X^{I_r}u \parallel_{S^{1,2}_{loc}(B(0,1))} &\leq&
C \Bigg( \parallel u
\parallel_{S^{1,2}_{loc}(B(0,2))} \\ &+& \parallel \tilde{f}
\parallel_{L^2(B(0,2))} + \parallel \tilde{f_i}
\parallel_{L^2(B(0,2))} \Bigg)
\end{eqnarray*}
where $\tilde{f}$, $\tilde{f_i}$ represent high order derivatives
on the original $f$, $f_i$.  The result follows once we apply
Corollary \ref{c1}.
\end{proof}

Finally, as a direct consequence to Corollary \ref{c1}, we have
that $u \in W^{k,2}_{loc}(B(0,1))$ for every k, where $W^{k,2}$
represents the usual Euclidean Sobolev space. Therefore, Corollary
\ref{MC} is immediate by the Sobolev Embedding Theorem.

\section{Sketch of Proof of Corollary $\ref{Gi}$}
As mentioned in the introduction, this is only a sketch of the
proof.  For further details we refer the reader to [CG] and [Gi].

\begin{proof}
The following corollary (Corollary $\ref{appcor}$) and inequality
that follows (see $\eqref{3.21}$ below) are a direct consequence
of the hypoellipticity result of the linear system in section 3
and will be used in the proof:
\begin{cor}\label{appcor}
Let G be a Carnot group, step r, and $\Omega \subset G$ an open
subset.  If $u \in S^{1,2}_{loc}(\Omega)$ is a weak solution of
the constant coefficient system
$$\sum_{i,j=1}^m \sum_{\beta=1}^N A_{ij}^{\alpha \beta} X_iX_j
u^\beta \,=\,0,\qquad \alpha =1,\ldots, N,$$ in $B(p_0,3R) \subset
\Omega$, then $u$ is smooth in $B(p_0,3R)$.  Moreover, there
exists a positive constant C such that $$\sup_{B(p_0,R)}(|u|^2 +
R^2 |Xu|^2 + R^4 \sum_{i,j=1}^m |X_iX_ju|^2) \leq C
\frac{1}{|B(p_0,2R)|}\int_{B(p_0,2R)}|u|^2 \,dp.$$
\end{cor}

Using Corollary $\ref{appcor}$, we can then prove the following
inequality holds for each $0 < r < R <2$, where C is a positive
constant and $G$, $\Omega$, and $u$ are the same as in Corollary
$\ref{appcor}$:
\begin{eqnarray}\label{3.21}
\int_{B(p_0,r)}|u-u_{(0,r)}|^2dp \leq
C(\frac{r}{R})^{Q+2}\int_{B(p_0,R)}|u-u_{0,R}|^2dp\,.
\end{eqnarray}

We use the following notation and prove the inequalities that
follow, assuming always that $u\in S_{loc}^{1,2}(\Omega)$ is a
weak solution to $\eqref{1.3}$. Set
\begin{eqnarray*}
U({p}_0,R)=\frac{1}{\vert{B({p}_0,R)\vert}}\int_{B({p}_0,R) \cap
\Omega}{\vert u(p)-u_{p_0,R}\vert}^2dp.
\end{eqnarray*}

The first step in the proof is to show that for each $M>0$ and
$0<\tau <1$, there exists $\epsilon_0$ and $R_0>0$ such that if
one has $\vert u_{p_0,R}\vert \leq M$ and $U(p_0,R)< \epsilon_0^2$
for $R\leq min(R_0,d(p_0,\partial \Omega))$ and for some $p_0 \in
\Omega$, then the following inequality holds:
\begin{eqnarray*}
U({p}_0,\tau R)\leq C \tau ^2 U({p}_0,R).
\end{eqnarray*}
The argument is by contradiction.  Set
$$
\upsilon ^n(q) = \epsilon _n^{-1}[u^n(p_n \delta _{R_n}(q)) -
u^n_{p_n,R_n}],
$$
so that
\begin{eqnarray}\label{4.8}
V^n( \epsilon ,1)= \frac{1}{|B_1|}\int_{B_1}|\upsilon ^n|^2 \,dq=1
\,,
\end{eqnarray} and assume
\begin{eqnarray}\label{4.9}
V^n(\epsilon,\tau) > 2C \tau^2\,.
\end{eqnarray}
Passing eventually to a subsequence and incorporating the
continuity assumptions on $A_{ij}^{\alpha \beta}$ along with the
hypothesis of the corollary, we obtain $$A_{ij}^{\alpha \beta}
(p_n \delta _{R_n}(q),\epsilon _n \upsilon ^n + u
^n_{p_n,R_n})\rightarrow B_{ij}^{\alpha \beta}\,.$$  Arguing as in
(4.10) pg. 25 in [CG], we then have that for every $\phi \in
C^\infty_0(B_1)$
\begin{eqnarray}
\sum_{\beta=1}^N \sum_{i,j=1}^m \int_{B_1} B_{ij}^{\alpha \beta}
X_i^i \upsilon^\beta X_i^1 \phi^\alpha \, dp \,=\,0\,,\quad \alpha
= 1,\ldots , N.
\end{eqnarray}
We can then apply inequality $\eqref{3.21}$ to get $$V(0,\tau)
\leq c \tau ^2 V(0,1)\,.$$  But following from $\eqref{4.8}$ and
$\eqref{4.9}$ we also must have that $$V(0,\tau) > 2C \tau ^2\,,$$
which is a contradiction.

An induction argument is used to show that, for every integer $k$,
\begin{eqnarray*}
U({p}_0,\tau^k R)\leq (2C \tau ^2)^k U({p}_0,R).
\end{eqnarray*}
It then follows from this inequality that
\begin{eqnarray*}
U(p,R)\leq CR^{2 \alpha},
\end{eqnarray*}
for each $R>0$ small enough and for $0<\alpha<1$ Consequently, we
now have that $u$ is H$\ddot{o}$lder continuous with exponent
$\alpha$ outside of a certain set that can be shown to have Haar
measure zero.
\end{proof}

\end{document}